\documentclass{conm-p-l}
\usepackage{array,epsfig,latexsym,tabularx}

\begin{document}


\newtheorem{example}{Example}[section]
\newtheorem{note}[example]{Note}
\newtheorem{theorem}[example]{Theorem}
\newtheorem{corollary}[example]{Corollary}
\newtheorem{definition}[example]{Definition}
\newtheorem{proposition}[example]{Proposition}
\newtheorem{algorithm}[example]{Algorithm}
\newtheorem{lemma}[example]{Lemma}
\newtheorem{problem}[example]{Problem}
\newtheorem{conjecture}[example]{Conjecture}


\newcommand{\bi}{\begin{itemize}}
\newcommand{\ei}{  \end{itemize}}
\newcommand{\be}{\begin{equation}}
\newcommand{\ee}{  \end{equation}}
\newcommand{\bea}{\begin{eqnarray}}
\newcommand{\eea}{  \end{eqnarray}} 
\newcommand{\leftgroup }{\left  \lgroup}
\newcommand{\rightgroup}{\right \rgroup}
\newcommand{\bw }{Boltzmann weight }
\newcommand{\bws}{Boltzmann weights }
\newcommand{\rai}{\rightarrow \, \infty}
\newcommand{\infinity}{\infty}
\newcommand{\plusmin}{\pm}
\newcommand{\ph}{\phantom}

 
\font\twelvesym=msbm10 at 12pt
\font\tensym=msbm10
\font\sevensym=msbm7
\font\fivesym=msbm5
\newfam\ssymfam
\textfont\ssymfam=\tensym
\scriptfont\ssymfam=\sevensym
\scriptscriptfont\ssymfam=\fivesym
\def\ssym{\fam\ssymfam\tensym}


\newcommand{\Z}{{\ssym Z}}
\newcommand{\N}{{\ssym N}}

\newcommand{\B}{{\mathcal B}}
\newcommand{\D}{{\mathcal D}}
\renewcommand{\H}{{\mathcal H}}
\newcommand{\M}{{\mathcal M}}
\renewcommand{\P}{{\mathcal P}}
\newcommand{\V}{{\mathcal V}}
\newcommand{\Sc}{{\mathcal S}}

\newcommand{\wt}{{\rm wt\,}}
\newcommand{\bwt}{wt \,}
\newcommand{\Proof}{\medskip\noindent {\it Proof: }}
\newcommand{\cqfd}{\hfill $\Box$ \medskip}
\newcommand{\boldm}{\mbox{\boldmath$m$}}
\newcommand{\boldn}{\mbox{\boldmath$n$}}
\newcommand{\bolde}{\mbox{\boldmath$e$}}
\newcommand{\boldC}{\mbox{\boldmath$C$}}
\newcommand{\sboldm}{\mbox{\boldmath$\scriptstyle m$}}
\newcommand{\sboldn}{\mbox{\boldmath$\scriptstyle n$}}
\newcommand{\sbolde}{\mbox{\boldmath$\scriptstyle e$}}
\newcommand{\sboldC}{\mbox{\boldmath$\scriptstyle C$}}
\newcommand{\boldlambda}{\mbox{\boldmath$\lambda$}}
\newcommand{\wombat}{\rule[-6pt]{0pt}{46pt}}






\def\makestrut#1#2{{\dimen12=#2
\divide\dimen12 by 4\dimen11=\dimen12\multiply\dimen11 by 3
\global\setbox#1=\hbox{\vrule height\dimen11 depth\dimen12 width0pt}}}

\newdimen\tadhdimen \newdimen\tabhdimen \newdimen\vdimen
\newdimen\smtadhdimen \newdimen\smtabhdimen
\newcount\splurt
\newbox\tadstrut \newbox\tabstrut
\newbox\smtadstrut \newbox\smtabstrut


\def\setyoungsize#1#2{          
  \tadhdimen=#1\tabhdimen=#1\advance\tabhdimen by -0.4truept%
  \vdimen=#2%
  \makestrut\tadstrut\vdimen
  \advance\vdimen by -0.4pt%
  \makestrut\tabstrut\vdimen}


\setyoungsize{18pt}{13pt}


\def\youngt#1{%
  \vcenter{\offinterlineskip
  \halign{&\copy\tadstrut\hbox to \tadhdimen{\hss$##$\hss}\cr #1}}}
\def\youngd#1{%
  \vcenter{\offinterlineskip
  \halign{&\vrule##&\copy\tabstrut\hbox to \tabhdimen{\hss$##$\hss}\cr #1}}}


\hyphenation{boson-ic 
             ferm-ion-ic 
	     para-ferm-ion-ic
             two-dim-ension-al
	     two-dim-ension-al}

\title{Path generating transforms}
\thanks{Research supported by the Australian Research Council (ARC)}
\dedicatory{Dedicated to Professor Richard Askey on the occasion 
            of his 65th birthday.}
\author{Omar~Foda}
\address{Department of Mathematics and Statistics,
             The University of Melbourne,
             Parkville, Victoria 3052, Australia.}
\email{foda@maths.mu.oz.au}
\author{Keith~S.~M.~Lee}
\address{Department of Mathematics and Statistics,
             The University of Melbourne,
             Parkville, Victoria 3052, Australia.}
\email{ksml@maths.mu.oz.au}
\author{Yaroslav~Pugai}
\address{Department of Mathematics and Statistics,
             The University of Melbourne,
             Parkville, Victoria 3052, Australia, \paragraph* 
             L.~D.~Landau Institute for Theoretical Physics,
             Russian Academy of Sciences, Moscow, Russia.}
\email{slava@itp.ac.ru}
\author{Trevor~A.~Welsh}
\address{Department of Mathematics and Statistics,
             The University of Melbourne,
             Parkville, Victoria 3052, Australia.}
\email{trevor@maths.mu.oz.au}
                 		      
\begin{abstract}
		       
We study combinatorial aspects of $q$-weighted, length-$L$ 
Forrester-Baxter paths, $\P^{p, p'}_{a, b, c}(L)$, where
$p, p', a, b, c\in\Z_{+}$, 
$0 < p < p'$, 
$0 < a, b, c < p'$, 
$c = b \pm 1$,
$L+a-b \equiv 0$ $(mod \, 2)$,
and 
$p$ and $p'$ are co-prime.

We obtain a bijection between $\P^{p, p'}_{a, b, c}(L)$ and 
partitions with certain prescribed hook differences.
Thereby, we obtain a new description of the $q$-weights of 
$\P^{p, p'}_{a, b, c}(L)$.
Using the new weights, and defining
$s_0$ and $r_0$ to be the smallest non-negative integers 
for which $|p s_0 - p' r_0|=1$,
we restrict the discussion to
$\P^{p, p'}_{s_0}$ $\equiv$ 
$\P^{p, p'}_{s_0, s_0, s_0 + 1}(L)$, 
and introduce two combinatorial transforms:
  
\begin{description}
\item[1] A Bailey-type transform $\B$: 
                    $\P^{p, p'}_{s_0}(L)$
	            $\rightarrow$
	            $\P^{p, p'+p}_{s_0 + r_0}(L')$, 
$L \leq L'$,

\item[2] A duality-type transform $\D$: 
                    $\P^{p,    p'}_{s_0}(L)$
		    $\rightarrow$
		    $\P^{p'-p, p'}_{s_0}(L)$.
\end{description}
We study the action of $\B$ and $\D$, as $q$-polynomial 
transforms on the $\P^{p, p'}_{s_0}(L)$ generating 
functions, $\chi^{p, p'}_{s_0}(L)$. In the limit 
$L \rightarrow \infinity$, $\chi^{p, p'}_{s_0}(L)$ reduces to 
the Virasoro characters, $\chi^{p, p'}_{r_0, s_0}$, of minimal 
conformal field theories $\M^{p, p'}$, or equivalently, 
to the one-point functions of regime-III Forrester-Baxter models. 
 
As an application of the $\B$ and $\D$ transforms, we re-derive 
the constant-sign expressions for $\chi^{p, p'}_{r_0, s_0}$,
first derived by Berkovich and McCoy.

\end{abstract}

\maketitle

\newpage

\setcounter{secnumdepth}{10}
\setcounter{section}{-1}

\section{Introduction}

\subsection{Motivation}

Many problems in exactly solvable models can be formulated as
problems in combinatorics. In particular, the computation of {\it
one-point functions} in two dimensional lattice models 
\cite{baxter-book}, can be reduced to the evaluation of the 
generating functions of certain combinatorial objects known as 
{\it one-dimensional configurations}. The purpose of this work 
is to study the combinatorics of a specific class of one-dimensional 
configurations called {\it paths}. In particular, we wish to 
study the paths that originate in computing one-point functions 
in an infinite series of lattice models introduced by Forrester 
and Baxter \cite{fb}\footnote{In fact, we restrict our attention 
to paths of {\it regime-III} Forrester-Baxter models. We refer 
the reader to \cite{fb} for definitions of the various regimes. 
{}From now on, we use the words \lq Forrester-Baxter models\rq, 
with the above restriction in mind.}. 

Roughly speaking, we would like to show how a set of paths 
that belongs to a certain model in an infinite series, can 
be systematically obtained from a simpler model in the same
series. The idea is to reduce the computation of one-point
functions, in any model, recursively, down to a computation
in the simplest possible model in the series. This latter
computation can be trivial.

The above idea is not new. It originates, at the level of $q$-series,
in the Bailey transform \cite{bailey}, and in its extensions by 
Andrews \cite{andrews-on-bailey}. Combinatorially, it also appears at 
the level of $q$-polynomials that count partition pairs \cite{burge}. 
Finally, and closest to the spirit of this work, it appears in 
\cite{ab,bressoud}, where certain infinite paths are $q$-counted. 
In this work, we formulate the above idea for certain finite length 
paths, and thereby we are able to obtain their generating functions. 

Our long term aim is to show that solutions in distinct models 
within a family, are not entirely independent, but are related 
to each other. Furthermore, if these solutions are understood 
combinatorially\footnote{For example, by showing that what one 
is doing is nothing but counting certain objects that satisfy 
certain properties.}, then one can obtain one from the other 
systematically, using combinatorial transforms. We anticipate 
that our approach will have applications to more general 
models, such as those based on the affine algebras $\hat{sl(n)}$.

\subsection{$\M^{p, p'}$ and $\chi^{p, p'}_{r, s}$}

Consider the minimal conformal field theories, $\M^{p, p'}$, of 
Belavin, Polyakov, and Zamolodchikov \cite{bpz}\footnote{For an 
introduction to conformal field theories, see \cite{dms-book}.}, 
where $p, p' \in \Z_{+}$ are co-prime and $1 < p < p'$.%
\footnote{In fact, when
we consider finite versions of these characters, our analysis 
naturally includes the cases $\{p=1,  p'\}$, which correspond to 
the {\it parafermionic} models \cite{dms-book}.}
$\M^{p, p'}$ are also the 
critical limits of the exactly solvable Forrester-Baxter models 
\cite{fb}\footnote{For an introduction to exactly solvable 
lattice models, see \cite{baxter-book}.}. 

The chiral half space of states $\M^{p, p'}$ \cite{bpz} may be 
considered as the union of a certain set of highest weight Virasoro 
modules of central charge $c = 1 - 6 (p - p')^2 / p p'$. Consider 
the characters, 
$\chi^{p, p'}_{r, s}$, where 
$p, p', r, s\in \Z_+$, 
$0 < r < p$, 
$0 < s < p'$, 
of the irreducible highest weight modules of the Virasoro algebra
which comprise $\M^{p, p'}$. In the language of conformal field 
theories, they are the simplest examples of {\it conformal blocks 
on the torus} \cite{dms-book}. In the language of lattice models, 
they are {\it one-point functions} up to a normalisation constant. 

Following Baxter's corner transfer matrix method \cite{baxter-book}, 
$\chi^{p, p'}_{r, s}$ can be interpreted as the generating functions of 
infinite length {\it paths}, $\P^{p, p'}_{r, s}$, that satisfy certain 
restrictions. These restrictions are encoded in the labels 
$p, p', r, s$ \cite{fb}. 
$\chi^{p, p'}_{r, s}$ can also be obtained as the generating functions 
of partitions with prescribed {\it hook differences} \cite{abbbfv}. 
There is no known physical or algebraic motivation for this latter
description.

\subsection{${\P}^{p, p'}_{a, b, c}(L)$}

The infinite length paths, ${\P}^{p, p'}_{r, s}$, can be regarded 
as the $L \rightarrow \infinity$ limit of length-$L$ paths 
${\P}^{p, p'}_{a, b, c}(L)$, where 
$0 < a, b, c < p'$, $c = b \pm 1$, $L+a-b \equiv 0$ $(mod \, 2)$. 
Here, $r,s$ are related to $a,b,c$ by
$r = \left\lfloor p c / p' \right\rfloor + (b - c + 1)/2$, 
$s = a$.  
In this paper, we work entirely at the level of
${\P}^{p, p'}_{a, b, c}(L)$, only taking the limit 
$L \rightarrow \infinity$ at the end.

\subsection{Purpose}

The purpose of this work is three-fold:

\begin{enumerate}

\item We describe a bijection between the finite length paths 
      $\P^{p, p'}_{a, b, c}(L)$, and partitions that satisfy 
      prescribed hook difference conditions \cite{abbbfv}. 
      Thereby, we obtain a new description of the $q$-weights 
      of $\P^{p, p'}_{a, b, c}(L)$. 

\item We introduce two combinatorial transforms, $\B$ and $\D$, 
      that act on finite length paths, and that can be used 
      recursively to generate 
      $\P^{p, p'}_{s_0}(L) \equiv \P^{p, p'}_{s_0, s_0, s_0 + 1}(L)$, 
      for all allowed $p, p'$, from 
      the combinatorially trivial $\P^{1, 3}_{1}(L')$.
      Here, $s_0$ is such that $s_0$ and $r_0$ are the smallest
      non-negative integers for which $|p s_0 - p' r_0|=1$, 
      and $L\ge0$ is even.

\item As an application, we compute the generating functions 
      of $\P^{p, p'}_{s_0}(L)$, and reproduce the constant sign 
      expressions for the characters $\chi^{p, p'}_{r_0, s_0}$, 
      that were first obtained in \cite{bm}.

\end{enumerate}

\subsection{Organisation}

Although the bijection between $\P^{p, p'}_{a, b, c}(L)$ and 
the corresponding partitions with prescribed hook differences 
is logically the starting point of this work, we relegate its
discussion and proof to an appendix.
This is so that we don't deviate from the main point of the 
paper, namely the combinatorics of the Forrester-Baxter paths 
and the $\B$ and $\D$ transforms.

In \S 1, we introduce the paths that we are interested in 
$q$-counting and their $q$-weights as originally defined
in \cite{fb}. We define extra structures on the paths, namely, 
the bands and their parities, an alternative prescription for
path weights that follows from the bijection described in 
Appendix A, the scoring and non-scoring vertices, and the 
striking sequence of a path. 

In \S 2, we introduce the $\B$-transform that maps
$\P^{p, p'}_{s_0}(L)$, 
into
$\P^{p, p'+ p}_{s_0 + r_0}(L')$
for various $L'$, where $s_0$ and $r_0$ are as defined above.  
We also define the particle content of a path.

In \S 3, we introduce the $\D$-transform that maps
$\P^{p, p'}_{s_0}(L)$
to
$\P^{p' - p, p'}_{s_0}(L)$. 

In \S 4, we digress to discuss the continued fraction description 
of the paths, and the related {\it zones}. We also discuss the 
the $\boldm\boldn$-system associated with $\P^{p, p'}_{s_0}(L)$,
and introduce a matrix that generalises the Cartan matrix of the
Lie algebra $A_t$, and that was first defined in \cite{bm}. We then 
describe what we mean by a sector
$\Sc(\hat{\boldn})\subset\P^{p, p'}_{s_0}(L)$, labelled
by $\hat{\boldn}\in\N^t$

In \S 5, we compute the constant-sign generating functions of 
$\P^{p, p'}_{s_0}(L)$. In particular, we re-derive the 
constant-sign expressions of the generating functions, first
obtained in \cite{bm}. Finally, we discuss the $L\rightarrow\infinity$
limit of these generating functions.

\section{Paths}

\subsection{Definitions}

Let $p$ and $p'$ be positive coprime integers for which $0<p<p'$.
Then, given $a,b,c,L\in\Z$ such that $1\le a,b,c\le p'-1$, $b=c\pm1$,
$L\ge0$, $L + a - b \equiv 0$ ($mod\,2$),
a path $h \in \P^{p, p'}_{a, b, c}(L)$ is a sequence
$h_0, h_1, h_2, \ldots, h_L,$ of integers such that:

\begin{enumerate}
\item $1 \le h_i \le p'-1$  for $0 \le i \le L$, 
\item $h_{i+1} = h_i \pm 1$ for $0 \le i <   L$, 
\item $a = h_0 , 
       b = h_L. $
\end{enumerate}

\noindent Note that the values of $p$ and $c$ do not feature in 
the above restrictions. As described below, they specify how the 
elements of $\P^{p, p'}_{a, b, c}(L)$ are weighted.

\subsection{Heights, segments, and vertices}

The integers $h_0,h_1,h_2,\ldots,h_L$ are readily depicted as a sequence 
of {\it heights} on a two-dimensional $L \times (p'-1)$ grid. Adjacent 
heights are connected by {\it line segments} passing from $(i,h_i)$ to 
$(i+1,h_{i+1})$ for $0 \le i < L$. 

Scanning the path from left to right, each of these line segments points 
either in the NE direction or in the SE direction. The following is a 
typical path in the set $\P^{3,11}_{5, 3, 4}(28)$:

\medskip
\begin{center}
\epsfig{file=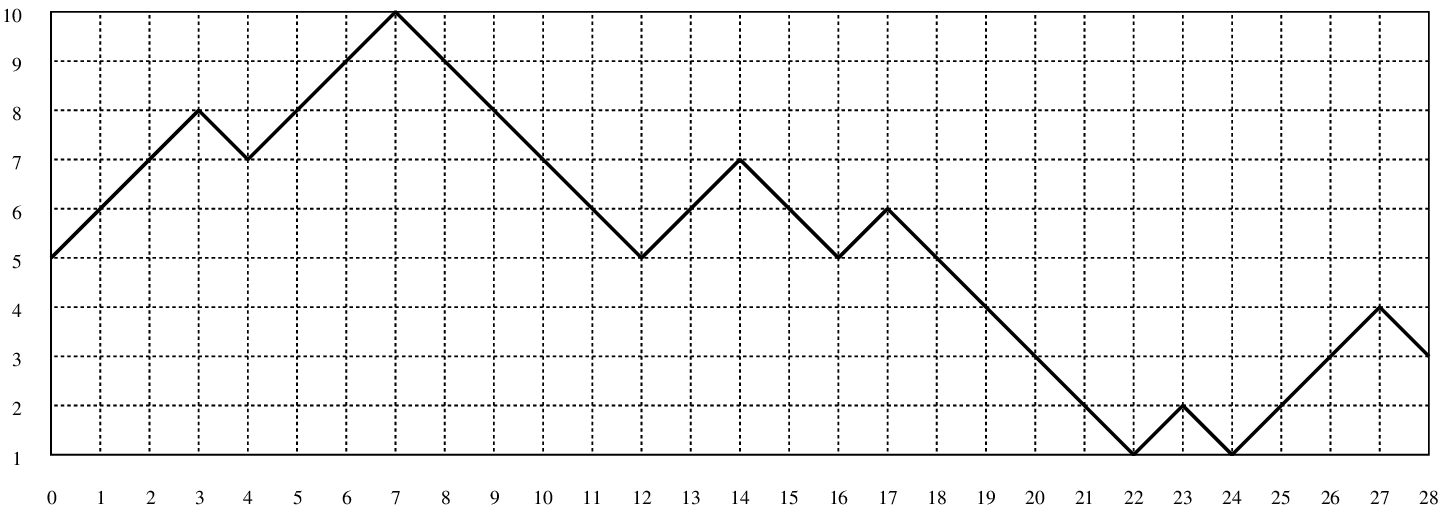, width=\linewidth}
\end{center}
\centerline{Figure\ 1.}
\medskip
\noindent

Two adjacent line segments, one passing 
from $(i-1,h_{i-1})$ 
to   $(i  ,h_{i  })$, 
and the other 
from $(i  ,h_{i  })$ 
to   $(i+1,h_{i+1})$, 
define a {\it vertex} $v_i$. 
There are four types of vertices. They appear as follows:

\begin{center}
\epsfig{file=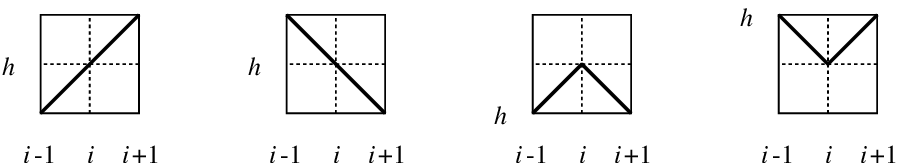}
\centerline{Figure\ 2.}
\end{center}

\noindent They will be referred to as a {\em straight-up} 
vertex, a {\em straight-down} vertex, a {\em peak-up} vertex 
and a {\em peak-down} vertex respectively.

\subsection{Forrester-Baxter weights}

In \cite{fb}, each vertex $v_i$ is assigned a weight $i c_{FB}$, 
where the local weight function $c_{FB}(h_{i-1},h_{i},h_{i+1})$ 
is defined by:

\begin{eqnarray*}
c_{FB}(h-1,h,h+1)&=&\phantom{+}1/2\,;\\
c_{FB}(h+1,h,h-1)&=&\phantom{+}1/2\,;\\
c_{FB}(h,h+1,h)&=&-\left\lfloor\frac{h(p'-p)}{p'}\right\rfloor\,;\\
c_{FB}(h,h-1,h)&=&\phantom{+}\left\lfloor\frac{h(p'-p)}{p'}\right\rfloor\,,
\end{eqnarray*}
\noindent where $\left\lfloor n \right\rfloor$ is the integer part of $n$.
A path $h$ is assigned a weight $\bwt(h)$ given by:

\begin{equation}\label{BwtDef}
\bwt(h)=\sum_{i=1}^L i c_{FB}(h_{i-1},h_{i},h_{i+1}),
\end{equation}
where we take $h_{L+1}=c$. In \cite{fb}, the generating function of 
$\P^{p, p'}_{a, b, c}(L)$ is defined to be\footnote{In \cite{fb}, 
this generating function is denoted by either $D_L(a,b,c)$ or
$D^{(k)}_L(a,b,c)$, where $k=\lfloor c(p'-p)/p'\rfloor$.}:

\begin{equation}\label{BgenDef}
\phi^{p, p'}_{a, b, c}(L) = 
\sum_{h \in \P^{p, p'}_{a, b, c}(L)} q^{\bwt(h)}.
\end{equation}

An expression for $\phi^{p,p'}_{a,b,c}(L)$ was derived 
in \cite[Theorem 2.3.1]{fb}. 
It turns out that there is a very convenient renormalisation
$\chi^{p,p'}_{a,b,c}(L)$ of $\phi^{p,p'}_{a,b,c}(L)$
(which we give explicitly in Appendix A), from which
the expression of \cite{fb} yields:

\begin{eqnarray}\label{FinRochaEq}
\chi^{p,p'}_{a,b,c}(L)&=&
\sum_{\lambda=-{\infinity}}^{\infinity}
q^{\lambda^{2} p p'+ \lambda (p'r-pa)}
\left[ {L \atop {{L+a-b} \over {2}}-p'\lambda} \right]_q\nonumber\\[0.5mm]
&&\qquad\qquad
-\sum_{\lambda=-\infinity}^\infinity
q^{(\lambda p+r)(\lambda p'+a)}
\left[ {L \atop {{L+a-b} \over {2}}-p'\lambda-a} \right]_q,
\end{eqnarray}

\noindent where 

\begin{equation}\label{groundstatelabel}
r=\lfloor pc/p'\rfloor+(b-c+1)/2,
\end{equation}

\noindent and, as usual, the Gaussian polynomial
$\left[ {A \atop B} \right]_q$ is defined to be:

\begin{equation}\label{gaussian}
\left[ {A \atop B} \right]_q = 
            \frac{\prod_{i=1}^{A  }(1 - q^i)}{\prod_{i=1}^{B  }(1 - q^i) 
                  \prod_{i=1}^{A-B}(1 - q^i)} 
\end{equation}
for $0\le B\le A$, and $\left[ {A \atop B} \right]_q=0$ otherwise.
In the limit $L\to\infinity$, we obtain

\begin{equation}\label{ChiLimEq}
\lim_{L\to\infinity}
\chi^{p,p'}_{a,b,c}(L) = \chi^{p,p'}_{r,a},
\end{equation}

\noindent where $r$ is defined in (\ref{groundstatelabel}) and

\begin{equation}\label{RochaEq}
\chi^{p, p'}_{r, s}=
{{1} \over {(q)_\infty}}\sum_{\lambda=-\infty}^\infty
(q^{\lambda^2pp'+\lambda(p'r-ps)}-q^{(\lambda p+r)(\lambda p'+s)})
\end{equation}

\noindent is, up to a normalisation, the Rocha-Caridi expression 
\cite{rocha} for the Virasoro character of central charge 
$c = 1 - {6(p' - p)^2}/{p p'}$
and conformal dimension $\Delta^{p, p'}_{r, s} = 
{\left( (p' r - p s)^2 - (p' - p)^2 \right)}/{4 p p'} $.
Therefore, $\chi^{p,p'}_{a,b,c}(L)$ provides a finite analogue of 
the character $\chi^{p,p'}_{r, a}$. 

The expression obtained above for $\chi^{p,p'}_{a,b,c}(L)$ is an 
alternating-sign $q$-polynomial. This expression is not combinatorial 
in the sense that we know that $\chi^{p,p'}_{a,b,c}(L)$ is a generating 
function, and therefore all its non-vanishing coefficients are positive. 
We shall seek constant-sign%
\footnote{For physical reasons, the alternating-sign 
expressions are also called {\it bosonic}. The constant-sign expressions 
are also called {\it fermionic}. The study of constant-sign expressions
for the Virasoro characters was initiated by the Stony Brook group.
For further details, and original references, we refer the reader to 
\cite{bm} and references therein.}
expressions for $\chi^{p,p'}_{a,b,c}(L)$, 
which in the limit $L\to\infinity$ will provide constant-sign
expressions for the Virasoro characters $\chi^{p, p'}_{r, a}$.

\subsection{Bands and parities}

In the path picture described above, there are $(p'-1)$ heights.
The regions between adjacent heights will be referred to as 
{\it bands}. There are $(p'-2)$ bands. We assign a {\it parity} 
to each band: a band that lies between heights $h$ and $(h+1)$ 
is {\em even} if
$\lfloor{hp/p'}\rfloor = \lfloor{(h+1)p/p'}\rfloor$, and {\em odd} 
otherwise. Scanning from below, the $r$th odd band lies between 
heights $\lfloor rp'/p\rfloor$ and $\lfloor rp'/p\rfloor+1$.
We will shade the odd bands more heavily than the even bands,
as shown in Fig.\ 3, where $p=3$ and $p'=11$.

\medskip
\begin{center}
\epsfig{file=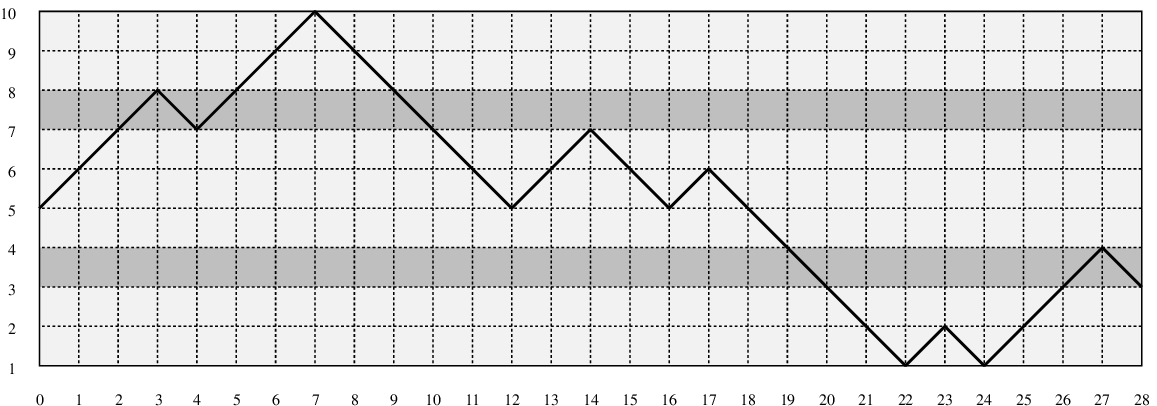, width=\linewidth}
\end{center}
\centerline{Figure\ 3.}
\medskip
\noindent

Since 
$\lfloor{hp/p'}\rfloor=0$ for $h=1$ and $\lfloor{hp/p'}\rfloor=p-1$
for $h=p'-1$, we deduce that there are $(p-1)$ odd bands and $(p'-p-1)$ 
even bands. Furthermore, we may also readily deduce that if $p'>2p$ then 
the odd bands occur is isolation, with an even band on both sides. If 
$p'<2p$, the reverse is true. Finally, it is easily seen that the odd/even 
band structure is invariant under an up/down reflection.

A parity may now be assigned to each vertex of a path: it is
the parity of the band in which the right edge lies.

\subsection{Alternative prescription for path weights}\label{PresSec}

The analysis of Appendix A shows that it is possible to assign
a weight $wt(h)$ to each path $h$ such that
$$
\chi^{p,p'}_{a,b,c}(L)
=\sum_{h\in{\P}^{p,p'}_{a,b,c}(L)} q^{wt(h)}.
$$
Here, we describe how $wt(h)$ may be simply calculated using the path
picture together with its band structure.
First, we define new coordinates on the picture as follows:
$$
x={{i-(h-a)} \over {2}},\qquad y={{i+(h-a)} \over {2}}.
$$
Thus, the $xy$-coordinate system has its origin at the path's initial
point, and is slanted at $45^o$ to the original $ih$-coordinate system.
Note that at each step in the path, either $x$ or $y$ is incremented
and the other is constant. In this system, the path depicted
in Fig.\ 3 has its first few coordinates at
$(0,0)$, $(0,1)$, $(0,2)$, $(0,3)$, $(1,3)$, $(1,4)$, $(1,5)$, $(1,6)$,
$(2,6)$, $\ldots$
Now, for the $i$th vertex, we define $c_i= c(h_{i-1},h_i,h_{i+1})$
according to the shape of the vertex and its parity.

\begin{table}[h]
\begin{center}
\begin{tabular}{|c|@{\hspace{3mm}}c@{\hspace{3mm}}|c|@{\hspace{3mm}}
c@{\hspace{3mm}}|}
\hline
Vertex&
${c}_i$&
Vertex&
${c}_i$\\
\hline\hline\wombat
{\epsfig{file=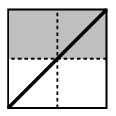}}&
\raisebox{14pt}[0pt]{x}&
{\epsfig{file=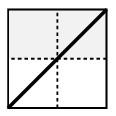}}&
\raisebox{14pt}[0pt]{0}\\
\hline\wombat
{\epsfig{file=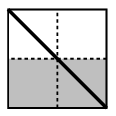}}&
\raisebox{14pt}[0pt]{y}&
{\epsfig{file=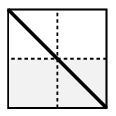}}&
\raisebox{14pt}[0pt]{0}\\
\hline\wombat
{\epsfig{file=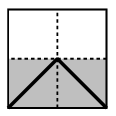}}&
\raisebox{14pt}[0pt]{0}&
{\epsfig{file=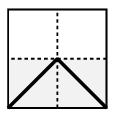}}&
\raisebox{14pt}[0pt]{x}\\
\hline\wombat
{\epsfig{file=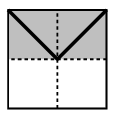}}&
\raisebox{14pt}[0pt]{0}&
{\epsfig{file=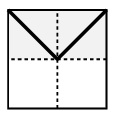}}&
\raisebox{14pt}[0pt]{y}\\
\hline
\end{tabular}
\end{center}
\end{table}

\noindent
Here the unshaded band can be either an even or an odd band (or in 
the lowermost four cases, not a band in the model at all). Note that 
each vertex shape only contributes in one parity case. We shall refer 
to those four vertices, with assigned parity, for which, in general, 
the contribution is non-zero, as {\em scoring} vertices. The other 
four vertices will be termed {\em non-scoring}.

We now define

\begin{equation}\label{WtDef}
wt(h)=\sum_{i=1}^L c_i.
\end{equation}

To illustrate this procedure, consider again the path $h$ depicted in
Fig.\ 3. The above table indicates that there are scoring vertices
at $i=2$, $7$, $9$, $12$, $14$, $16$, $17$, $19$,
$22$, $23$, $24$, $26$.
This leads to
$$
wt(h)=0+1+6+6+6+8+8+9+9+13+10+14
=90.
$$

\subsection{Striking sequence of a path}

Consider each path $h$ as a sequence of straight lines, alternating
in direction between NE and SE. Reading from the left, let the lengths
of these lines be 
$w_1$, $w_2$, $w_3,\ldots,w_l,$ for some $l$, so that each $w_i>0$ 
and $w_1+w_2+\cdots+w_l=L(h)$, where $L(h)$ is the length of $h$.
For each of these lines, the last vertex will be considered to be
part of the line but the first will not. Then, the $i$th of these 
lines contains $w_i$ vertices, the first $w_i-1$ of which are
straight vertices. Then write $w_i=a_i+b_i$ so that $b_i$ is the 
number of scoring vertices in the $i$th line. The striking sequence 
of $h$ is then the array:

\begin{displaymath}\label{HseqDef}
\left(\begin{array}{ccccc}
  a_1&a_2&a_3&\cdots&a_l\\ b_1&b_2&b_3&\cdots&b_l
 \end{array}\right).
\end{displaymath}
We define $m(h)=\sum_{i=1}^l a_i$, whereupon $\sum_{i=1}^l b_i=L(h)-m(h)$.
We also define $\beta(h)=(b_1+b_3+\cdots)-(b_2+b_4+\cdots)$.

For example, for the path shown in Fig.\ 3, the striking sequence is:

$$
\def\qua{\hskip5pt}
\left(
{2\qua1\qua2\qua3\qua1\qua1\qua0\qua3\qua0\qua0\qua2\qua1\atop
 1\qua0\qua1\qua2\qua1\qua1\qua1\qua2\qua1\qua1\qua1\qua0}
\right).
$$
In this case, $m(h)=16$, and $\beta(h)=0$.

\begin{lemma}\label{WtHashLem}
Let the path $h$ have the striking sequence
$\left({a_1 \atop b_1}\:{a_2 \atop b_2}\:{a_3 \atop b_3}\:
 {\cdots\atop\cdots}\:{a_l\atop b_l} \right),$
with $w_i=a_i+b_i$ for $1\le i\le l$.
Then
$$
wt(h)=\sum_{i=1}^l b_i(w_{i-1}+w_{i-3}+\cdots+w_{1+i\bmod2}).
$$
\end{lemma}

\Proof First assume that the first $w_1$ segments of $h$ are in 
the NE direction. Then, for $i$ odd, the $i$th line is in the NE 
direction and its $x$-coordinate is $w_2+w_4+\cdots+w_{i-1}$. By 
the prescription of the previous section, and the definition of 
$b_i$, this line thus contributes $b_i(w_2+w_4+\cdots+w_{i-1})$ 
to the weight $wt(h)$ of $h$. Similarly, for $i$ even, the $i$th 
line is in the SE direction and contributes 
$b_i(w_1+w_3+\cdots+w_{i-1})$ to $wt(h)$.
This proves the lemma if the first segments are in the NE direction.
The reasoning is almost identical for the other case.
\cqfd

\subsection{Restricting the endpoints of paths}

In the rest of this work, we restrict our attention to the set of
paths $\P^{p, p'}_{s_0}(L) \equiv \P^{p, p'}_{s_0, s_0, s_0+1}(L)$,
where $L\ge0$ is even and $s_0$ is defined to be such that $s_0$ and 
$r_0$ are the smallest non-negative integers for which 
$|p s_0 - p' r_0| = 1$. In the limit $L \rightarrow \infinity$, the 
corresponding generating functions reduce to the Virasoro characters 
related to the models $\M^{p, p'}$, with smallest possible highest 
weights. The line $h=s_0$ in the path picture will be referred to as
the {\em ground-line}.

Using the fact that $p$ and $p'$ are co-prime, it is straightforward
to deduce that there do exist $s_0$ and $r_0$ satisfying
$0\le r_0\le s_0\le p'-1$ such that $|p s_0 - p' r_0| = 1$.
In fact, in the case $p=1$ we immediately obtain $s_0=1$ and $r_0=0$,
and in the case $p=p'-1$ we immediately obtain $s_0=1$ and $r_0=1$.
Otherwise, if $1<p<p'-1$, then necessarily $1<s_0<p'-1$.
Moreover, if $ps_0-p'r_0=1$, so that $ps_0/p'=r_0+1/p'$, the $h=s_0$
line is immediately below the $r_0$th odd band, and above
an even band.
On the other hand, if $ps_0-p'r_0=-1$, so that $ps_0/p'=r_0-1/p'$,
the $h=s_0$ line is immediately above the $r_0$th odd band, and below
an even band.
We make use of this information to derive the following technical
result for later convenience.

\begin{lemma}\label{EasyBetaLem}
For all $h^*\in{\P}^{p,p'}_{s_0}(L)$, we have
$\beta(h^*)=0$.
\end{lemma}
\Proof We first define a {\em flip} transformation of a path.
This consists of, exchanging two consecutive segments that form
a down-peak for two that form an up-peak, or vice-versa.
Thus a flip appears as follows:

\medskip
\begin{center}
\epsfig{file=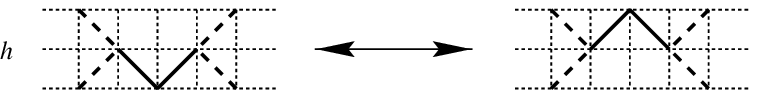}
\end{center}
\medskip

Note that the segment preceding the two that are changed
(which appears only if the peak is not the first vertex)
and the segment succeeding the two that are changed,
may each be either up or down, and each band may be even or odd.
Consideration of all sixteen cases shows that if the two paths
$h^\prime$ and $h^{\prime\prime}$ differ by a flip transformation,
then $\beta(h^\prime)=\beta(h^{\prime\prime})$.
In the case where the peak is the first vertex of the path,
this holds if and only if the two bands are of opposite parity.
This is always so for the $h=s_0$ case that we are dealing with
here.
(Note that if $h=1$, this type of flip transformation is not
valid and is not used.)

Now, consider the unique path $h^{(0)}$ for which the
first $L$ segments all lie in the $s_0$th band.
We immediately obtain $\beta(h^{(0)})=0$.
Now, as is easily seen, any path
$h^*\in{\P}^{p,p'}_{s_0}(L)$
can be obtained from $h^{(0)}$ by a sequence of flip
transformations.
Therefore $\beta(h^*)=\beta(h^{(0)})=0$, as required.
\cqfd
\medskip

\section{The $\B$-transform}\label{BTranSec}

In this section, we introduce the {\it $\B$-transform}\footnote{Our 
$\B$-transform is basically a generalisation of a transform 
introduced by Agarwal and Bressoud in \cite{ab,bressoud}. This, 
in turn, was motivated by the works of Bailey \cite{bailey} 
and Burge \cite{burge}. The transform in
\cite{ab,bressoud} acts on the infinite paths that pertain to 
the cases where $p=2$ and $p'=2k+1$ for $k \in \Z_{+}$. In this 
work, we extend the transform of \cite{ab,bressoud}, so that it 
acts on finite paths. We also generalise it to all co-prime 
$p, p'$. However, whereas the analysis of \cite{ab,bressoud} 
applies to $\P^{2, 2k+1}_{a, k, k+1}$, for all $a$, we restrict 
consideration to a single value of $a$ in the current paper.}:
the first of the two combinatorial transforms studied in this 
work. 

The $\B$-transform comprises three separate steps which we refer 
to as the $\B_1$, $\B_2$ and $\B_3$-transforms. $\B_1$ maps 
$\P^{p, p'}_{s_0}(L)$ (injectively) to 
$\bigcup_{L^{(0)}}\P^{p, p'+p}_{s_0 + r_0}(L^{(0)})$.
$\B_2$ lengthens paths in a simple manner, in particular mapping 
$\P^{p, p'}_{s_0}(L)$ to $\P^{p, p'}_{s_0}(L+2k)$ for 
$k\ge0$. $\B_3$ deforms, in a particular manner, a path within
$\P^{p, p'}_{s_0}(L)$. As we will see, taking the paths generated
by the combined action of $\B_1$ followed by $\B_2$ as input
to $\B_3$, the paths generated by $\B_3$ are naturally indexed 
by certain partitions $\lambda$ (See Appendix \ref{PartSec} for 
a definition of a partition).

The $\B$-transform comprising $\B_1$ followed by $\B_2$ followed
by $\B_3$, and involving the two parameters $k$ and $\lambda$,
will be denoted $\B(k,\lambda)$, and
thus maps from $\bigcup_L\P^{p, p'}_{s_0}(L)$ to
$\bigcup_{L'}\P^{p, p'+p}_{s_0+r_0}(L')$.
As we will show, this map: $(h,k,\lambda)\mapsto h'$ is
actually a bijection.

\subsection{The $\B_1$-transform}

The definition of the $\B_1$-transform involves the band structure
of $\P^{p, p'}_{s_0}(L)$. First note that the band structure of 
$\P^{p, p'+p}_{s_0 + r_0}(L')$ is easily obtained from that of 
$\P^{p, p'}_{s_0}(L)$. In the two cases, the number of odd bands 
is the same.  Since the $r$th odd band for $\P^{p, p'}_{s_0}(L)$ 
has its lower edge at height $\lfloor rp'/p\rfloor$ and that of 
$\P^{p, p'+p}_{s_0 + r_0}(L')$ has its lower edge at 
$\lfloor r(p'+p)/p\rfloor=\lfloor rp'/p\rfloor+r$,
we see that the distance between the odd bands has increased by
exactly one, with the height of the lowermost having also increased
by one. Note that, if $p>1$, the starting point is on the upper or 
lower edge of the $r_0$th odd band both before and after the 
$\B_1$-transform. If $p=1$, the starting point remains at $h=1$.

The image of the path is now obtained by examining the sequence of 
vertex types and inserting an extra one immediately prior to each 
scoring vertex. An extra straight-up vertex is inserted immediately
prior to each odd straight-up vertex and each even peak-down vertex,
and an extra straight-down vertex is inserted immediately prior to 
each odd straight-down vertex and each even peak-up vertex. In view 
of the odd bands having separated by one unit, and the change of the 
starting point, we see that the shapes and parities of the scoring 
vertices are naturally preserved under this transform.

For example, the following $\P^{3,8}_{3}(16)$ path (here $L=16$ and
we show $h_{L+1}=c=4$):

\medskip
\begin{center}
\epsfig{file=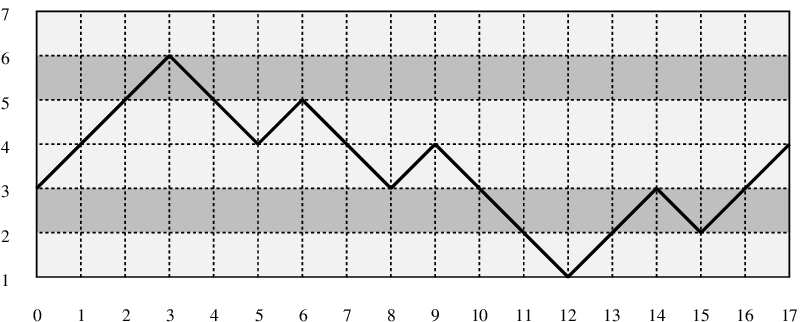}
\end{center}
\medskip

\noindent
maps to this $\P^{3,11}_{4}(24)$ path:

\medskip
\begin{center}
\epsfig{file=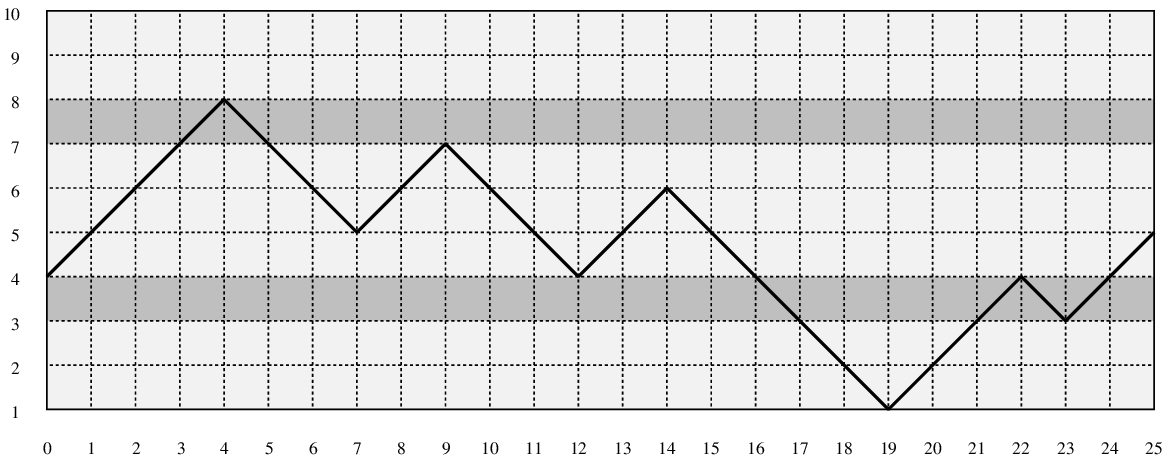}
\end{center}
\medskip

\noindent
Note that the path obtained from a $\B_1$-transform is such that there
are no two consecutive scoring vertices, and that the first vertex is
non-scoring.

\begin{lemma}\label{BHashLem}

Let $h \in \P^{p, p'}_{s_0}(L)$
have striking sequence
$\left({a_1\atop b_1}\:{a_2\atop b_2}\:{a_3\atop b_3}\:
 {\cdots\atop\cdots}\:{a_l\atop b_l} \right),$
and let
$h^{(0)} \in \P^{p,p'+p}_{s_0 + r_0}(L^{(0)})$ be obtained from
the action of the $\B_1$-transform on $h$.
Then $h^{(0)}$ has striking sequence:
\begin{displaymath}
\left(\begin{array}{ccccc}
  a_1+b_1&a_2+b_2&a_3+b_3&\cdots&a_l+b_l\\ b_1&b_2&b_3&\cdots&b_l
 \end{array}\right),
\end{displaymath}
$m(h^{(0)})=L$ and $L^{(0)}=2L-m(h)$.
\end{lemma}

\Proof This follows directly from the definition of the
striking sequences and the action of the $\B_1$-transform.
\cqfd
\medskip

\begin{lemma}\label{BWtLem}
Let $h \in{\P}^{p,p'}_{s_0}(L)$ 
and
$h^{(0)} \in \P^{p,p'+p}_{s_0 + r_0}(L^{(0)})$ be the path obtained 
by the action of the $\B_1$-transform on $h$.
Then
$$
wt(h^{(0)})=wt(h)+{ {1} \over {4}}(L-m)^2,
$$
where $m=m(h)$.
\end{lemma}

\Proof Let $h$ have striking sequence
$\left({a_1\atop b_1}\:{a_2\atop b_2}\:{a_3\atop b_3}\:
 {\cdots\atop\cdots}\:{a_l\atop b_l} \right),$
Then Lemmas \ref{BHashLem} and \ref{WtHashLem} show that
\begin{eqnarray*}
wt(h^{(0)})-wt(h)&=&(b_1+b_3+b_5+\cdots)(b_2+b_4+b_6+\cdots)\\
&=&{{1} \over {4}}((L-m)^2-\beta(h)^2),
\end{eqnarray*}
the second equality resulting because
$L-m=b_1+b_2+\cdots+b_l$ and
$\beta(h)=(b_1+b_3+b_5+\cdots)-(b_2+b_4+b_6+\cdots)$.
The lemma now follows on using Lemma \ref{EasyBetaLem}.
\cqfd
\medskip

\subsection{The $\B_2$-transform}

Let $p'>2p$ so that $h^{(0)}\in\P^{p,p'}_{s_0}(L')$ is a path for which
there are no two neighbouring odd bands.  By {\em inserting a particle} 
into $h^{(0)}$, we mean displacing $h^{(0)}$ two squares to 
the right and inserting two even edges. Since the path starts at 
height $s_0$ and the line $h=s_0$ borders both an even and an odd 
band (if $s_0>1$), this is possible in exactly one way. In this way, 
we obtain a path $h^{(1)}$ of length $L'+2$.
By repeating this process we may obtain 
a path $h^{(k)}\in\P^{p,p'}_{s_0}(L'+2k)$ by successively inserting $k$ 
particles into $h^{(0)}$. We say that $h^{(k)}$ has been obtained 
by the action of a $\B_2$-transform on $h^{(0)}$. When we need to 
show the dependence on $k$ explicitly, we refer to it as 
a $\B_2(k)$-transform.

\begin{lemma}\label{WtShiftLem}
Let $h \in \P^{p,p'}_{s_0}(L)$.
Apply a $\B_1$-transform to $h$ to obtain the path
$h^{(0)} \in \P^{p,p'+p}_{s_0 + r_0}(L^{(0)})$. 
Then obtain $h^{(k)}\in\P^{p,p'+p}_{s_0 + r_0}(L^{(k)})$ by
applying a $\B_2(k)$-transform to $h^{(0)}$.
If $m^{(k)}=m(h^{(k)})$,
then $L^{(k)}=L^{(0)}+2k$, $m^{(k)}=m^{(0)}$ and
$$
wt(h^{(k)})=wt(h)+{ {1} \over {4}}(L^{(k)}-m^{(k)})^2.
$$
\end{lemma}
\Proof
That $L^{(k)}=L^{(0)}+2k$ follows immediately from the
definition of a $\B_2$-transform.
Lemma \ref{BWtLem} yields:
\begin{eqnarray*}
wt(h^{(0)})&=&wt(h)+{{1}\over {4}}(L-m(h))^2\\
&=&wt(h)+{{1} \over {4}}\left(L^{(0)}-m(h^{(0)})\right)^2,
\end{eqnarray*}
the second equality following from Lemma \ref{BHashLem}.
Let the striking sequence of $h^{(0)}$ be
$\left({a_1\atop b_1}\:{a_2\atop b_2}\:
 {\cdots\atop\cdots}\:{a_l\atop b_l} \right),$
whereupon that of $h^{(1)}$ is
$\left({0\atop1}\:{0\atop1}\:{a_1\atop b_1}\:{a_2\atop b_2}\:
 {\cdots\atop\cdots}\:{a_l\atop b_l} \right),$
Then, $m(h^{(1)})=m(h^{(0)})$ and Lemma \ref{WtHashLem} shows that
$wt(h^{(1)})=wt(h^{(0)})+L^{(0)}-m(h^{(0)})+1$.
Repeated application then yields
$m(h^{(k)})=m(h^{(0)})$ and
$$
wt(h^{(k)})=wt(h^{(0)})+k\left(L^{(0)}-m(h^{(0)})\right)+k^2.
$$
Then, using $L^{(k)}=L^{(0)}+2k$, the lemma follows.
\cqfd
\medskip

\subsection{Particle moves}

In this section, we restrict attention to those $\P^{p,p'}_{s_0}(L)$ 
for which $p'>2p$, and specify six types of local deformations of 
a path. These deformations will be known as {\em moves}.\footnote{In 
the special cases $p = 1, 2$, these moves are equivalent to those 
considered in \cite{ab,bressoud,warnaar}. In the special cases 
$p'=p+1$, our results are equivalent to similar results obtained in 
\cite{berkovich,warnaar,df}.}  
In each of the six cases, 
a particular sequence of four segments of a path is changed to a 
different sequence, the remainder of the path being unchanged. The 
moves are as follows --- the path portion to the left of the arrow is 
changed to that on the right:

\bigskip
\centerline{\epsfbox{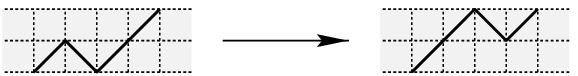}}
\nobreak
\centerline{Move.\ 1.}
\medskip
\centerline{\epsfbox{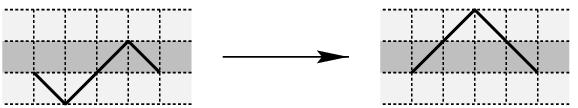}}
\nobreak
\centerline{Move.\ 2.}
\medskip
\centerline{\epsfbox{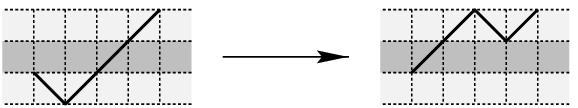}}
\nobreak
\centerline{Move.\ 3.}
\bigskip
\centerline{\epsfbox{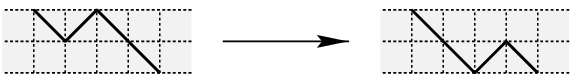}}
\nobreak
\centerline{Move.\ 4.}
\bigskip
\centerline{\epsfbox{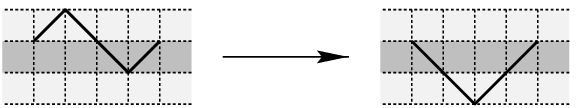}}
\nobreak
\centerline{Move.\ 5.}
\bigskip
\centerline{\epsfbox{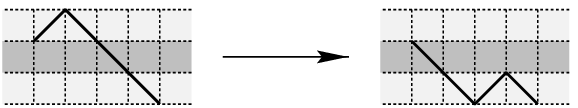}}
\nobreak
\centerline{Move.\ 6.}
\bigskip

\noindent Since $p'>2p$, each odd band is straddled by a pair
of even bands. Thus, there is no impediment to enacting moves
2 and 5 for paths $\P^{p,p'}_{s_0}(L)$.

Note that moves 4--6 are inversions of moves 1--3. Also note 
that moves 2 and 3 (likewise moves 5 and 6) may be considered 
to be the same move since in the two cases, the same sequence 
of three edges is changed.

\begin{lemma}\label{WtChngLem}
Let $h$ be a path for which four consecutive segments are as in 
one of diagrams on the left above. Let $\hat h$ be that obtained 
from $h$ by changing those segments according to the move. Then
$$
wt(\hat h)=wt(h)+1.
$$
Additionally, $L(\hat h)=L(h)$ and $m(\hat h)=m(h)$.
\end{lemma}
\Proof For each case, take the $xy$-coordinate of the leftmost 
point of this portion of a path to be $(x_0,y_0)$. Now consider 
the contribution to the weight of the three vertices in question 
before and after the move (although the vertex immediately before 
those considered may change, its contribution doesn't). In each 
of the six cases, the contribution is $x_0+y_0+1$ before the move 
and $x_0+y_0+2$ afterwards. Thus the first statement holds. The final 
statement is obtained by inspecting all six moves.
\cqfd
\medskip

Now observe that for each of the moves specified above, the sequence 
of path segments before the move consists of an adjacent pair of 
scoring vertices followed by a non-scoring vertex. When $p'>2p$ so 
that there are no two adjacent odd bands, all such combinations of 
vertices are present amongst the six moves. The specified move then 
consists of replacing such a combination with a non-scoring vertex 
followed by two scoring vertices. It is useful to interpret this as 
the pair of adjacent scoring vertices having moved rightward by one 
step. As anticipated above, we refer to a pair of adjacent scoring
vertices as a particle. Thus each of the six moves here is a particle
moving to the right by one step.

\subsection{The $\B_3$-transform}

Since each of the moves described above moves a pair of scoring
vertices to the right by one step, we see that a succession of such 
moves is possible until the pair is followed by another scoring 
vertex. If this itself is followed by yet another scoring vertex,
we forbid further movement. However, if it is followed by a non-scoring 
vertex, further movement is allowed after considering the latter two 
of the three consecutive scoring vertices to be the particle (instead 
of the first two).

As above, let $h^{(k)}$ be a path resulting from a $\B_2$-transform
inserting $k$ particles into a path that itself is the image of
a $\B_1$ transform. We now consider moving these $k$ particles.

\begin{lemma}\label{GaussLem} 
There is a bijection between the paths obtained by moving the
particles in $h^{(k)}$ and the partitions $\lambda$ with at
most $k$ parts, none of which exceeds $m=m(h^{(k)})$.
This bijection is such that if $h$ is the bijective image of 
a particular $\lambda$ then
$$
wt(h)=wt(h^{(k)})+wt(\lambda),
$$
where $wt(\lambda)=\lambda_1+\lambda_2+\cdots+\lambda_k$.
Additionally, $L(h)=L(h^{(k)})$ and $m(h)=m(h^{(k)})$.
\end{lemma}

\Proof
Since each particle moves by traversing a non-scoring vertex,
and there are $m$ of these to the right of the rightmost particle 
in $h^{(k)}$, and there are no consecutive scoring vertices to its 
right, this particle can make $\lambda_1$ moves to the right, with 
$0\le\lambda_1\le m$. Similarly, the next rightmost particle can 
make $\lambda_2$ moves to the right with $0\le\lambda_2\le\lambda_1$.
Here, the upper restriction arises because the two scoring vertices 
would then be adjacent to those of the first particle. Continuing in 
this way, we obtain that all possible final positions of the particles 
are indexed by $\lambda=(\lambda_1,\lambda_2,\ldots,\lambda_k)$ with 
$m\ge\lambda_1\ge\lambda_2\ge\cdots\ge\lambda_k\ge0$, that is, by 
partitions of at most $k$ parts with no part exceeding $m$. Moreover, 
since by Lemma \ref{WtChngLem} the weight increases by one for each 
move, the weight increase after the sequence of moves specified by 
a particular $\lambda$ is equal to $wt(\lambda)$.
The final statement also follows from Lemma \ref{WtChngLem}.
\cqfd
\medskip

\noindent We say that a path obtained by moving the particles
in $h^{(k)}$ has been obtained by the action of a $\B_3$-transform.
If we wish to show the specific dependence on $\lambda$, we refer 
to the transform as a $\B_3(\lambda)$-transform.

Having defined $\B_1$, $\B_2(k)$ for $k\ge0$ and $\B_3(\lambda)$
for $\lambda$ a partition with at most $k$ parts, we now define a 
$\B(k,\lambda)$-transform as the composition 
$\B(k,\lambda)=\B_3(\lambda)\circ\B_2(k)\circ\B_1$.

\subsection{Particle content of a path}

In this section let $p'>2p$ so that there are no two adjacent odd bands.  
Let $h'\in\P^{p,p'}_{s_0}(L')$.

\begin{lemma}\label{UniqueLem} There is a unique triple $(h,k,\lambda)$
where $h\in\P^{p,p'-p}_{s_0-r_0}(L)$ for some $L$, such that
the action of a $\B(k,\lambda)$-transform on $h$ results in $h'$.
\end{lemma}

\Proof This is proved by reversing the constructions described
in the previous sections. Locate the leftmost pair of consecutive 
scoring vertices in $h'$, and move them leftward by reversing the 
particle moves, until they occupy the first two positions. Now 
ignoring these two vertices, do the same with the next leftmost 
pair of consecutive scoring vertices, moving them leftward until 
they occupy the third and fourth positions. Continue in this way 
until all consecutive scoring vertices occupy the leftmost positions 
of the path. Say there are $2k$ of them (that this number is even is
implied by the starting point being between an odd and an even band),
and denote this path $h^{(k)}$. Clearly $h'$ results from $h^{(k)}$ 
by a $\B_3(\lambda)$-transform for a particular $\lambda$ with at most 
$k$ parts.

Removing the first $2k$ segments of $h^{(k)}$ produces a path
$h^{(0)}\in\P^{p,p'}_{s_0}(L'-2k)$ which has no two consecutive
scoring vertices. Moreover, $h^{(k)}$ arises by the action of
a $\B_2(k)$-transform on $h^{(0)}$.

Finally, since $h^{(0)}$ has by construction no pair of
consecutive scoring vertices, and none at the first vertex,
we may remove a non-scoring vertex before every scoring
vertex to obtain a path $h\in\P^{p,p'-p}_{s_0-r_0}(L)$ for
some $L$, from which $h^{(0)}$ arises by the action of a
$\B_1$-transform. The lemma is then proved.
\cqfd
\medskip

\noindent
The value of $k$ obtained above will be referred to as the
particle content of $h'$.

\section{The $\D$-transform}

The {\it $\D$-transform}\footnote{A similar duality-type transform 
appears in \cite{andrews-on-bailey,bm}.} is defined to act on each 
$h \in \P^{p,p'}_{s_0}(L)$ 
to yield a path in $\P^{p'-p,p'}_{s_0}(L)$. It is easily seen that
the band structure of 
$\P^{p'-p,p'}_{s_0}(L)$ is obtained from that of $\P^{p,p'}_{s_0}(L)$ 
simply by replacing odd bands by even bands and vice-versa. The action 
on $h \in \P^{p,p'}_{s_0}(L)$ yields a path 
$\hat h \in \P^{p'-p,p'}_{s_0}(L)$ with exactly the same sequence of
integer heights, i.e., $\hat h_i=h_i$ for $0\le i\le L$.

\begin{lemma}\label{DHashLem}
Let 
$h \in \P^{p,p'}_{s_0}(L)$ have striking sequence
$\left({a_1\atop b_1}\:{a_2\atop b_2}\:{a_3\atop b_3}\:
 {\cdots\atop\cdots}\:{a_l\atop b_l} \right),$
and let $\hat h \in \P^{p'-p,p'}_{s_0}(L)$ be obtained from
the action of the $\D$-transform on $h$.
Then $\hat h$ has striking sequence:
\begin{displaymath}
\left(\begin{array}{ccccc}
  b_1&b_2&b_3&\cdots&b_l\\ a_1&a_2&a_3&\cdots&a_l
 \end{array}\right).
\end{displaymath}
Moreover, $m(\hat h)=L(h)-m(h)$ and $L(\hat h)=L(h)$.
\end{lemma}

\Proof This follows directly from the definition of the
striking sequences after noting that the action of the $\D$-transform
exchanges odd bands for even bands and vice-versa.
\cqfd
\medskip

\begin{lemma}\label{DWtLem}
Let $h\in{\P}^{p,p'}_{s_0}(L)$, and $\hat h\in \P^{p'-p,p'}_{s_0}(L)$
be obtained by the action of a $\D$-transform on $h$. Then
$$
wt(h)+wt(\hat h)={{1}\over {4}}L^2.
$$
\end{lemma}

\Proof
Let $h$ have striking sequence
$\left({a_1\atop b_1}\:{a_2\atop b_2}\:{a_3\atop b_3}\:
 {\cdots\atop\cdots}\:{a_l\atop b_l} \right),$
and let $w_i=a_i+b_i$ for $1\le i\le l$.
Then, using Lemmas \ref{WtHashLem} and \ref{DHashLem}, we obtain
\begin{eqnarray*}
wt(h)+wt(\hat h)&=& 
\sum_{i=1}^l b_i(w_{i-1}+w_{i-3}+\cdots+w_{1+i\bmod2})\\
&&\qquad
+ \sum_{i=1}^l a_i(w_{i-1}+w_{i-3}+\cdots+w_{1+i\bmod2})\\[0.5mm]
&=&
\sum_{i=1}^l w_i(w_{i-1}+w_{i-3}+\cdots+w_{1+i\bmod2})\\[0.5mm]
&=&
(w_1+w_3+w_5+\cdots)(w_2+w_4+w_6+\cdots).
\end{eqnarray*}
The lemma then follows because
$(w_1+w_3+w_5+\cdots)+(w_2+w_4+w_6+\cdots)=L$ and
$(w_1+w_3+w_5+\cdots)-(w_2+w_4+w_6+\cdots)=0$ because
the start and endpoints of the paths have equal heights.
\cqfd
\medskip

To obtain the particle content of a path
$h'\in\P^{p,p'}_{s_0}(L')$, where $p'<2p$,
we first perform a $\D$-transform on the
path $h'$ to obtain a path $\hat h'\in \P^{p'-p,p'}_{s_0}(L')$.
The particle content of $h'$
is defined to be equal to that of $\hat h'$.

In what follows, we obtain the model $\P^{p,p'}_{s_0}(L)$ by using 
a certain sequence of $\B$- and $\D$-transforms. This sequence is 
determined by the continued fraction expansion of $p'/p$.

\section{Continued fractions and the $\boldm\boldn$-system}

\subsection{Basic definitions}

If $p'$ and $p$ are positive co-prime integers and
$$
{p'\over p}=
{c_0+{\displaystyle\strut 1 \over \displaystyle c_1 +
{\displaystyle\strut 1 \over \displaystyle c_2 +
{\displaystyle\strut 1 \over {\lower-5pt\hbox{$\vdots$}
\over \displaystyle\strut c_{n-1} +
{\displaystyle\strut 1 \over \displaystyle c_n}}}}}}
$$
with $c_0\ge0$, $c_i\ge1$ for $0<i<n$, and $c_n\ge2$,
then $(c_0,c_1,c_2,\ldots,c_n)$ is said to be the
{\it continued fraction} for $p'/p$.\footnote{Although 
we make no use of this fact, it may be easily shown
(\cite[Lemma 3.9]{flw}) that with $s_0$ and $r_0$ defined 
to be the smallest non-negative co-prime integers such 
that $|p s_0 - p' r_0|=1$, then $s_0/r_0$ has continued 
fraction $(c_0,c_1,\ldots,c_{n-1})$ if $n>0$.}

For later convenience, we now determine how the continued fractions
of $p'/p$ change under the $\B$- and $\D$-transforms.

\begin{lemma}\label{CFTranLem}

For positive co-prime integers $p'$ and $p$, let $p'/p$ have
continued fraction $(c_0,c_1,c_2,\ldots,c_n)$.

\begin{enumerate}

\item The continued fraction of $(p'+p)/p'$ is
$(c_0+1,c_1,c_2,\ldots,c_n)$. 

\item If $p'>2p$ then the continued fraction of $p'/(p'-p)$ is
$(1,c_0-1,c_1,c_2,\ldots,c_n)$.
\end{enumerate}

\end{lemma}

\Proof The first result follows immediately. So does the second
after writing $p'/(p'-p)=1+1/(p'/p-1)$.
\cqfd
\medskip

\subsection{Zones}\label{ZoneSec}

For what follows, it is convenient to partition indices into {\it zones}.
To this end, given $p$ and $p'$ with $p'/p$ having continued fraction 
$(c_0,c_1,\ldots,c_n)$, define\footnote{The $t_{n+1}$ defined here differs 
from that defined in \cite{bm}.}
\begin{equation}\label{ZoneEq}
t_\mu=-1-\delta_{\mu,n+1}+\sum_{i=0}^{\mu-1} c_i\qquad(\mu=0,1,\cdots,n+1).
\end{equation}
Then we say that the index $j$ with $0\le j\le t_{n+1}$ is in zone $\mu$
if $t_{\mu}<j\le t_{\mu+1}$. Thus, there are $n+1$ zones. Note that for 
$0\le\mu<n$, zone $\mu$ contains $c_\mu$ indices, and zone $n$ contains 
$c_n-1$ indices.

We define $t=t_{n+1}$ and refer to it as the {\it rank} of $p'/p$. We also 
define the rank and number of zones of ${\P}^{p,p'}_{a,b,c}(L)$ to be equal 
to the rank and number of zones of $p'/p$ respectively. Then Lemma 
\ref{CFTranLem} shows that the rank is incremented under a $\B$-transform 
while the number of zones is constant, and the rank is constant under a 
$\D$-transform while the number of zones is incremented.

Let $p'>2p$ so that $t_1>0$. We now see that if we begin with 
$\P^{1,3}_{1}(L')$ and apply a sequence of $t-1$ $\B$-transforms,
with, for $1\le\mu\le n$, the $(t-t_\mu)$th $\B$-transform immediately 
preceded by a $\D$-transform, then we obtain elements of 
$\bigcup_{L}\P^{p,p'}_{s_0}(L)$.\footnote{In fact, it is admissible 
to extend this sequence and consider it as passing from the trivial
$\P^{1,2}_{1}(0)$ to $\bigcup_{L}\P^{p,p'}_{s_0}(L)$.} In what 
follows, we will use this sequence of transforms to obtain the paths 
that we wish to enumerate.

\subsection{The $\boldm\boldn$-system}

For each pair of positive co-prime integers $p,p'$, we now define
the associated {\it $\boldm\boldn$-system}%
\footnote{The $\boldm\boldn$-systems defined here were first defined
in \cite{bm}, where the interrelationship between $\boldm$ and $\boldn$
was derived by an altogether different method.}.
Let $p'/p$ have rank $t$.
The $\boldm\boldn$-system is then a set of $t$ linear equations
defining an interdependence between two $t$-dimensional vectors%
\footnote{In this paper, the vectors $\boldm$, $\boldn$ and $\hat{\boldn}$
should be considered as column vectors. However, for typographical
convenience, we shall express their components in row vector form.}
$\hat{\boldn}=(n_1,n_2,\ldots,n_{t})$
and $\boldm=(m_0,m_1,\ldots,m_{t-1})$.
The equations are given by, for $1\le j\le t$:
\begin{eqnarray}
&&m_{j-1}-m_{j+1}=m_{j}+2n_{j}\quad
\hbox{if $j=t_\mu,\quad\mu=1,2,\ldots,n$;}
\label{MNEq1}\\
&&m_{j-1}+m_{j+1}=2m_j+2n_j\quad\hbox{otherwise,}
\label{MNEq2}
\end{eqnarray}
where we set $m_{t}=m_{t+1}=0$.

Note that if each $n_i$ is a non-negative integer then each $m_j$ is an 
even non-negative integer. Also note that it is possible to eliminate 
$m_j$ for $1\le j\le t$ from the above equations to obtain:
$$
\sum_{i=1}^{t} l_in_i={m_0\over 2},
$$
for certain positive integer values $l_1,l_2,\ldots, l_{t}$. These 
values are referred to as {\it string lengths} in \cite{bm} and feature 
prominently in the analysis there. We don't require them in the present 
paper.  In what follows, $m_0$ will be identified with the length $L$ 
of a path.

As an example, consider the case where $p=9$ and $p'=31$.
Here, the continued fraction of $31/9$ is
$(3,2,4)$, whereupon $n=2$, $t_1=2$, $t_2=4$ and $t=t_3=7$.
The $\boldm\boldn$-system of equations yields the following:
\begin{eqnarray*}
m_7&=&0;\\
m_6&=&2n_7;\\
m_5&=&2n_6+4n_7;\\
m_4&=&2n_5+4n_6+6n_7;\\
m_3&=&2n_4+2n_5+6n_6+10n_7;\\
m_2&=&2n_3+4n_4+2n_5+8n_6+14n_7;\\
m_1&=&2n_2+2n_3+6n_4+4n_5+14n_6+24n_7;\\
L\:=\:m_0&=&2n_1+4n_2+2n_3+8n_4+6n_5+20n_6+34n_7.
\end{eqnarray*}

In addition to the above $t$ equations (not counting $m_t=0$),
it will be useful to
set $m_{-1}=0$ and to use the appropriate (\ref{MNEq1}) or (\ref{MNEq2})
in the case $j=0$ to obtain one further equation.
Thence define the $t$-dimensional vector
$\boldn=(n_0,n_1,\ldots,n_{t-1})$.

The significance of the $\boldm\boldn$-system will be revealed after we
examine sequences of $\B$- and $\D$-transforms, and insert a number of 
particles at each stage.

\subsection{The generalised Cartan matrix}

The above expressions for the $\boldm\boldn$-system may be
conveniently expressed in matrix form.
The matrix so involved is a generalisation of the Cartan matrix of
the Lie algebra of type A.

Given $p$ and $p'$, define $t_\mu$ as in (\ref{ZoneEq}) and let 
$t=t_{n+1}$. Now let $\boldC$ be the $t\times t$ tri-diagonal 
matrix with entries $\boldC_{ij}$ for $0\le i,j\le t-1$ where, 
when the indices are in this range,

\begin{displaymath}
\begin{array}{cccl}
\boldC_{j,j-1}=-1,
&\boldC_{j,j}=1,
&\boldC_{j,j+1}=\phantom{-}1,
&\hbox{if $j=t_\mu,\quad\mu=1,2,\ldots,n$;}\\
\boldC_{j,j-1}=-1,
&\boldC_{j,j}=2,
&\boldC_{j,j+1}=-1,
&\hbox{$0\le j<t$ otherwise.}
\end{array}
\end{displaymath}

\begin{lemma}\label{CarLem}
For fixed $p$ and $p'$, let $\boldm$ and $\boldn$
satisfy the $\boldm\boldn$-system.
Then:
$$
2\boldn=-\boldC\boldm.
$$
\end{lemma}

\Proof This follows immediately from the definition of $\boldC$ and
equations (\ref{MNEq1}) and (\ref{MNEq2}) for $0\le j<t$.
\cqfd
\medskip

\begin{corollary}\label{CarCor}
For fixed $p$ and $p'$, let $\boldm$ and $\boldn$
satisfy the $\boldm\boldn$-system.
Then, on setting $L=m_0$:
\begin{displaymath}
\sum_{i=1}^{t} m_in_i=
\left\{ \begin{array}{l@{\rule[3pt]{0pt}{10pt}\quad:\quad}l}
 -{1\over 2}\boldm^T\boldC\boldm + {1\over 2} Lm_1 + {1\over 2} L^2
  &\mbox{if }t_1=0;\\
 \strut -{1\over 2}\boldm^T\boldC\boldm - {1\over 2} Lm_1 + L^2
  &\mbox{if }t_1>0.
 \end{array} \right.
\end{displaymath}
\end{corollary}

\Proof Using Lemma \ref{CarLem},
$\sum_{i=1}^{t} m_in_i=
-{1\over 2}\boldm^T\boldC\boldm-m_0n_0$.
In the case $t_1=0$, expression (\ref{MNEq1}) gives $2n_0=-m_1-m_0$.
In the case $t_1>0$, expression (\ref{MNEq2}) gives $2n_0=m_1-2m_0$.
The result then follows after substituting $m_0=L$.
\cqfd
\medskip

\subsection{Sectors} In this section, we use the sequence of $\B$-
and $\D$-transforms detailed at the end of Section \ref{ZoneSec},
to pass from $\P^{1,3}_{1}(L)$ to $\P^{p,p'}_{s_0}(L')$.
First we prove two short lemmas.

\begin{lemma}\label{MNsys1Lem}
Let $h\in\P^{p,p'}_{s_0}(L)$ and let $h'\in\P^{p,p'+p}_{s_0+r_0}(L')$
arise through the action of a $\B(n',\lambda)$-transform on $h$.
If $m=m(h)$ and $m^\prime=m(h^\prime)$, then $m^\prime=L$ and
$$
L^\prime+m=2m^\prime+2n^\prime.
$$
\end{lemma}

\Proof By Lemma \ref{BHashLem}, the action of the $\B_1$-transform
yields a path $h^{(0)}$ of length $2L-m$, with $m(h^{(0)})=L$.
Thereupon, using Lemma \ref{WtShiftLem}, the action of a
$\B_2(n')$-transform on $h^{(0)}$ results in a path $h^{(k)}$
of length $2L-m+2n'$ for which $m(h^{(k)})=L$.
Since a $\B_3(\lambda)$-transform does not change either
property, the lemma follows.
\cqfd
\medskip

\begin{lemma}\label{MNsys2Lem}
Let $h\in\P^{p,p'}_{s_0}(L)$ and let $h'\in\P^{p'-p,2p'-p}_{2s_0-r_0}(L')$
arise through the action of a $\D$-transform on $h$ followed by
a $\B(n',\lambda)$-transform.
If $m=m(h)$ and $m^\prime=m(h^\prime)$, then $m^\prime=L$ and
$$
L^\prime-m=m^\prime+2n^\prime.
$$
\end{lemma}

\Proof Let $\hat h\in\P^{p'-p,p'}_{s_0}(L)$ result from the action
of the $\D$-transform on $h$.
By Lemma \ref{DHashLem}, $m(\hat h)=L-m$ and $L(\hat h)=L$.
Then, using Lemma \ref{MNsys1Lem}, $m^\prime=L$ and
$L^\prime+(L-m)=2m^\prime+2n^\prime$.
The required expressions now follow.
\cqfd
\medskip
 
We now define subsets $\Sc^{p,p'}(\hat{\boldn})$
of $\P^{p,p'}_{s_0}(L)$
indexed by $\hat{\boldn}=(n_1,n_2,\ldots,n_t)\in\N^t$, which we
refer to as sectors%
\footnote{We take $0\in\N$.}.
For $p'>2p$, the sector $\Sc^{p,p'}(\hat{\boldn})$ indexed by
$\hat{\boldn}$ is obtained as follows.
First let $\boldm=(m_0,m_1,\ldots,m_{t-1})$ be obtained from
the $\boldm\boldn$-system of $p'/p$ and set $m_t=0$.
Now let
$\boldlambda=(\lambda^{(1)},\lambda^{(2)},\ldots,\lambda^{(t-1)})$
be a sequence of partitions such that for each $1\le i<t$,
$\lambda^{(i)}$ has at most $n_i$ parts, none of
which exceeds $m_i$.

We now create a sequence $(h^{(t)},h^{(t-1)},\ldots,h^{(1)})$ of
paths as follows.
First let $h^{(t)}\in\P^{1,3}_1(2n_t)$.
Note that this path is unique.
Now for $1\le j<t$, obtain the path $h^{(j)}$ from
$h^{(j+1)}$ as follows.
If $j\ne t_\mu$ for all $\mu=1,2,\ldots,n$, obtain $h^{(j)}$
by applying a $\B(n_j,\lambda^{(j)})$-transform to $h^{(j+1)}$.
In the case where $j=t_\mu$, for $\mu=1,2,\ldots,n$,
apply a $\D$-transform to $h^{(j+1)}$ immediately before applying
a $\B(n_j,\lambda^{(j)})$-transform.

\begin{lemma}\label{SectorLem} Let $p'>2p$ and let the sequence
$(h^{(t)},h^{(t-1)},\ldots,h^{(1)})$ of paths be obtained as above.
Then, for each $1\le j\le t$, we have
$m(h^{(j)})=m_j$ and $L(h^{(j)})=m_{j-1}$.
\end{lemma}

\Proof The result follows immediately in the case $j=t$
since $L(h^{(t)})=2n_t=m_{t-1}$ and $m(h^{(t)})=0$ by
direct inspection of $h^{(t)}$.
Now proceed by (decreasing) induction on $t$.
Let $1\le k<t$ and assume that the result holds for $j=k+1$,
so that $m(h^{(k+1)})=m_{k+1}$ and $L(h^{(k+1)})=m_k$.
If $k\ne t_\mu$ for all $\mu=1,2,\ldots,n$,
Lemma \ref{MNsys1Lem} immediately shows that
$m(h^{(k)})=L(h^{(k+1)})=m_k$ and $L(h^{(k)})=2m_k+2n_k-m_{k+1}$,
which, by (\ref{MNEq2}), is $m_{k-1}$.
The case $k=t_\mu$ for $\mu=1,2,\ldots,n$, follows
similarly from Lemma \ref{MNsys2Lem} and (\ref{MNEq1}).
The required result therefore holds in the case $j=k$, whereupon
the lemma follows by induction.
\cqfd
\medskip

For $p'>2p$, we now define the sector $\Sc^{p,p'}(\hat{\boldn})$
to be the set of all paths obtained as above for all $\boldlambda$
within the specified constraints.
Lemma \ref{SectorLem} then shows that for $h\in\Sc^{p,p'}(\hat{\boldn})$,
the values $L(h)$ and $m(h)$ are constant and given by $m_0$ and
$m_1$ respectively.
In particular, $\Sc^{p,p'}(\hat{\boldn})\subset\P^{p,p'}_{s_0}(m_0)$.
For $p'<2p$, the sector $\Sc^{p,p'}(\hat{\boldn})$ is defined to be
precisely the image of the sector
$\Sc^{p'-p,p'}(\hat{\boldn})\subset\P^{p'-p,p'}_{s_0}(m_0)$
under the $\D$-transform.

Now for $p'>2p$, given $h\in\P^{p,p'}_{s_0}(L)$ we can reverse
the above construction (if $p'<2p$ we simply apply a $\D$-transform
first). Let $h^{(1)}=h$ and define a sequence 
$(h^{(t)},h^{(t-1)},\ldots,h^{(1)})$ of paths,
$\hat{\boldn}=(n_1,n_2,\ldots,n_t)$ and a sequence
$\boldlambda=(\lambda^{(1)},\lambda^{(2)},\ldots,\lambda^{(t-1)})$
of partitions as follows.
For $1\le j< t$, obtain $h^{(j+1)}$ from $h^{(j)}$ as follows.
Lemma \ref{UniqueLem} shows that there is a unique triple
$(h',k,\lambda)$ for which $h^{(j)}$ arises from the action of the
$\B(k,\lambda)$-transform on $h'$.
Set $n_j=k$ and $\lambda^{(j)}=\lambda$.
If $j\ne t_\mu$ for $\mu=1,2,\ldots,n$, set $h^{(j+1)}=h'$.
Otherwise, set $h^{(j+1)}$ to be the result of the $\D$-transform
acting on $h'$.
In this way, we find that $h$ is present in one and only one
sector, namely $h\in\Sc^{p,p'}(\hat{\boldn})$.
Moreover, the sequence $\boldlambda$ of partitions is uniquely
determined by $h$.

\section{Generating functions and character formulae}

\subsection{Constant-sign generating functions for paths}

In this section, we combine the techniques of the previous sections
to calculate constant-sign generating functions for $\P^{p,p'}_{s_0}(L)$.
The first step is to determine such a function for all paths in
a sector. Given co-prime $p,p'$, let $p'/p$ have rank $t$ and
let $\hat{\boldn}=(n_1,n_2,\ldots,n_t)\in\N^t$.
Now define the sector generating function to be:
\begin{equation}
S^{p,p'}(\hat{\boldn})=\sum_{h\in\Sc^{p,p'}(\hat{\sboldn})} q^{\wt(h)}.
\end{equation}

\begin{lemma}\label{StepLem} For $p'>2p$, let
$\hat{\boldn}=(n_1,n_2,\ldots,n_t)$ and
$\boldm=(L,m_1,\ldots,m_{t-1})$ satisfy the
$\boldm\boldn$-system of $p'/p$.
Then, if we set $\hat{\boldn}'=(n_2,n_3,\ldots,n_t)$,
$$
S^{p,p'}(\hat{\boldn})=
q^{{1\over4}(L-m_1)^2}
\left[{m_1+n_1\atop n_1}\right]_q
S^{p,p'-p}(\hat{\boldn}').
$$
\end{lemma}

\Proof Lemma \ref{SectorLem} shows that for each
$h\in\Sc^{p,p'}(\hat{\boldn})$, $L(h)=L$ and $m(h)=m_1$.
Then, Lemma \ref{UniqueLem} shows
that there is a bijection between $\Sc^{p,p'}(\hat{\boldn})$
and pairs $(h',\lambda)$ where
$h'\in\Sc^{p,p'-p}(\hat{\boldn}')$ and $\lambda$ is a partition
with at most $n_1$ parts, none of which exceeds $m_1$.
This bijection is such that
$h$ results from the action of the $\B(n_1,\lambda)$-transform
on $h'$.
Moreover, Lemmas \ref{WtShiftLem} and \ref{GaussLem} show that
$\wt(h)=\wt(h')+{1\over4}(L-m_1)^2+\wt(\lambda)$.
Therefore, if we let $\langle n_1,m_1\rangle$ denote the set of
partitions that have at most $n_1$ parts,
none of which exceeds $m_1$,
\begin{eqnarray*}
S^{p,p'}(\hat{\boldn})&=&
\sum_{h'\in\Sc^{p,p'-p}(\hat{\sboldn}')}
\sum_{\lambda\in\langle n_1,m_1\rangle}
q^{{1\over4}(L-m_1)^2+\wt(\lambda)+\wt(h')}\\
&=&
q^{{1\over4}(L-m_1)^2}
\left(\sum_{\lambda\in\langle n_1,m_1\rangle} q^{\wt(\lambda)}\right)
\left(\sum_{h'\in\Sc^{p,p'-p}(\hat{\sboldn}')} q^{\wt(h')}\right).
\end{eqnarray*}
Since the generating function for all partitions $\lambda$ with at
most $n_1$ parts, none of which exceeds $m_1$,
is $\left[{m_1+n_1\atop n_1}\right]_q$
(see \cite{andrews-red-book} for example),
the lemma now follows.
\cqfd
\medskip

\begin{lemma}\label{GenLem}
Let $p'$ and $p$ be positive co-prime integers and let $t$ be
the rank of $p'/p$.
Let $\hat{\boldn}=(n_1,n_2,\ldots,n_t)$
and $\boldm=(L,m_1,\ldots,m_{t-1})$ satisfy the
$\boldm\boldn$-system of $p'/p$.
If $p'>2p$, then
\begin{equation}\label{SecEq1}
S^{p,p'}(\hat{\boldn})=
q^{ {1\over 4}L(L-m_1) -{1\over 2}\sum_{j=1}^{t} m_jn_j}
\prod_{j=1}^{t-1}
\left[{m_j+n_j\atop n_j}\right]_q;
\end{equation}
and if $p'<2p$, then
\begin{equation}\label{SecEq2}
S^{p,p'}(\hat{\boldn})=
q^{ {1\over 4}Lm_1 -{1\over 2}\sum_{j=1}^{t} m_jn_j}
\prod_{j=1}^{t-1}
\left[{m_j+n_j\atop n_j}\right]_q.
\end{equation}
\end{lemma}

\Proof We first prove that, for a given $t$, expression (\ref{SecEq2})
follows from expression (\ref{SecEq1}). So let $p'<2p$, with $p'/p$
having rank $t$ and assume that (\ref{SecEq1}) holds for rank $t$.
In particular, it holds for the sector
$\Sc^{p'-p,p'}(n_1,n_2,\ldots,n_{t})$.
Notice that if $p'<2p$ then the $\boldm\boldn$-system is the same
for $p'/(p'-p)$ as it is for $p'/p$.
Using
$$
\left[{m+n}\atop n\right]_q
=q^{mn}\left[{m+n}\atop n\right]_{q^{-1}},
$$
expression (\ref{SecEq1}) may be rewritten
$$
S^{p'-p,p'}(\hat{\boldn})
=
q^{ {1\over 4}L(L-m_1) +{1\over 2}\sum_{j=1}^{t} m_jn_j}
\prod_{j=1}^{t-1}
\left[{m_j+n_j\atop n_j}\right]_{q^{-1}}.
$$
Since the action of a ${\D}$-transform on
$\Sc^{p'-p,p'}(\hat{\boldn})$ yields
exactly the set $\Sc^{p,p'}(\hat{\boldn})$,
Lemma \ref{DWtLem} gives:
$$
S^{p,p'}(\hat{\boldn})
=
q^{{1\over 4}L^2-\left({1\over 4}L(L-m_1)
               +{1\over 2}\sum_{j=1}^{t} m_jn_j\right)}
\prod_{j=1}^{t-1}
\left[{m_j+n_j\atop n_j}\right]_q,
$$
whereupon (\ref{SecEq2}) follows.

We now prove (\ref{SecEq1}) for $p'>2p$. 
In the case $p=1$ and $p'=3$, there is a unique path
of length $L$, which comprises exactly $n_1=L/2$ particles.
Its weight is easily found to be $n_1^2$.
Therefore $S^{1,3}(\hat{\boldn})=q^{{1\over4}L^2}$.
Since $p'/p=3/1$ has continued fraction $(3)$,
then in this case $t=1$ and thence $m_1=0$.
Therefore, (\ref{SecEq1}) holds for rank $t=1$.

We now proceed by induction on the rank of $p'/p$.
Assume that (\ref{SecEq1}) and (\ref{SecEq2}) hold for all
$p'/p$ of rank $k-1$, where $k>1$.
Now let $p'/p$ have rank $k$ so that if
$p'/p$ has continued fraction $(c_0,c_1,\ldots,c_n)$
then $c_0+c_1+\ldots+c_n-2=k$.
Note that $(p'-p)/p$ has continued fraction
$(c_0-1,c_1,\ldots,c_n)$ and hence rank $k-1$.
Let $\hat{\boldn}'=(n_2,n_3,\ldots,n_t)$.
We consider separately the two cases $c_0>2$ and $c_0=2$ for which
$t_1>1$ and $t_1=1$ respectively.
In the $c_0>2$ case, the induction hypothesis gives:
$$
S^{p,p'-p}(\hat{\boldn}')=
q^{ {1\over 4}m_1(m_1-m_2) -{1\over 2}\sum_{j=2}^{k} m_jn_j}
\prod_{j=2}^{k-1}
\left[{m_j+n_j\atop n_j}\right]_q.
$$
Then Lemma \ref{StepLem} yields:
$$
S^{p,p'}(\hat{\boldn})=
q^{ {1\over 4}(L-m_1)^2 }
\left[{m_1+n_1\atop n_1}\right]_q
q^{ {1\over 4}m_1(m_1-m_2) -{1\over 2}\sum_{j=2}^{k} m_jn_j}
\prod_{j=2}^{k-1}
\left[{m_j+n_j\atop n_j}\right]_q.
$$
Expression (\ref{SecEq1}) now follows for the case $t=k$,
after noting that
$$
(L-m_1)^2=L(L-m_1)-m_1(L-m_1)=L(L-m_1)-m_1(m_1-m_2+2n_1),
$$
where the second equality arises from (\ref{MNEq2}) with $j=1$.

In the $c_0=2$ case, $(p'-p)/p$ has continued fraction
$(1,c_1,c_2,\ldots,c_n)$ and thus $2p<(p'-p)$.
Then combining (\ref{SecEq2}) at $t=k-1$ with 
Corollary \ref{StepLem}, results in:
$$
S^{p,p'}(\hat{\boldn})=
q^{ {1\over 4}(L-m_1)^2 }
\left[{m_1+n_1\atop n_1}\right]_q
q^{ {1\over 4}m_1m_2 -{1\over 2}\sum_{j=2}^{k} m_jn_j}
\prod_{j=2}^{k-1}
\left[{m_j+n_j\atop n_j}\right]_q.
$$
In this case, expression (\ref{SecEq1}) now follows for the case $t=k$,
after noting that
$$
(L-m_1)^2=L(L-m_1)-m_1(L-m_1)=L(L-m_1)-m_1(m_2+2n_1),
$$
where here, the second equality arises from (\ref{MNEq1}) with $j=t_1=1$.

Thus, expressions (\ref{SecEq1}) and (\ref{SecEq2}) at $t=k$ follow
from those at $t=k-1$. The lemma then follows by induction.
\cqfd
\medskip

We are now able to give a constant-sign generating function for 
all paths in ${\P}^{p,p'}_{s_0,s_0,s_0+1}(L)$ and thus provide 
a constant-sign expression for the finitised character 
$\chi^{p,p'}_{r_0,s_0}(L)$.\footnote{We note that this expression 
combines the two expressions for $p' < 2p$, and $p' > 2p$ given 
in \cite{bm}.}

\begin{theorem}\label{LFinThrm1}
Let $p$ and $p'$ be positive co-prime integers and $r_0$ and $s_0$
the smallest positive integers such that $\vert ps_0-p'r_0\vert=1$.
Then for even $L\ge0$:
\begin{eqnarray*}
\chi^{p,p'}_{s_0,s_0,s_0+1}(L)
&=&
\sum_{\sboldm}
q^{{1\over 4}\sboldm^T\sboldC\sboldm-{1\over 4}L^2}
\prod_{j=1}^{t-1}
\left[{m_j-{1\over 2}(\boldC\boldm)_j\atop m_j}\right]_q,
\end{eqnarray*}
where the summation is over vectors
$\boldm=(L,m_1,m_2,\ldots,m_{t-1})$
with $m_j\in2\Z$ and $m_j\ge0$ for $j=1,2,\ldots,t-1$.
\end{theorem}

\Proof Since the two cases $p'>2p$ and $p'<2p$ correspond to
$t_1>1$ and $t_1=0$ respectively, Corollary \ref{CarCor} shows that
the exponents in (\ref{SecEq1}) and (\ref{SecEq2}) are each
equal to ${1\over 4}\boldm^T\boldC\boldm-{1\over 4}L^2$.
The current theorem therefore follows after summing over
all sectors $(n_1,n_2,\ldots,n_t)$ which give $m_0=L$, and writing
$n_j=-{1\over 2}(\boldC\boldm)_j$.
That each $m_j$ is even, follows inductively from (\ref{MNEq1})
and (\ref{MNEq2}).
\cqfd
\medskip

\subsection{Character formulae}

In Theorem \ref{LFinThrm1}, we note that no explicit restrictions
are given for the $m_j$ (other than being even and non-negative).
In fact, there do exist implicit restrictions,
that arise because for
$$
\left[{m_j+n_j\atop m_j}\right]_q
$$
($n_j=-{1\over 2}(\boldC\boldm)_j$) to be non-zero requires $n_j\ge0$.
In order to more conveniently impose these restrictions and then to take
the limit $L\to\infinity$, we now change to a further set of variables.
This set will also enable a comparison with the results of \cite{flw}
to be made.
The variables will be $\lambda^{(\mu)}_i$ for $0\le\mu\le n$
and $1\le i\le c_\mu-2\delta_{\mu,n}$.
It will turn out that we need only consider the cases for
which each $\lambda^{(\mu)}$ is actually a partition.
It also turns out that we need only consider those cases for which
certain $\lambda^{(\mu)}_1$ and certain $wt(\lambda^{(\mu)})$ are
bounded above by values that depend on $L$ and the partitions
$\lambda^{(0)}, \lambda^{(1)},\ldots,\lambda^{(\mu-1)}$.
These two bounds, if imposed, will be denoted $\lambda^{(\mu)}_0$
and $w_\mu$ respectively.

For the purpose of having uniform expressions, set
$(e_0,e_1,\ldots,e_{n-1},e_n)=(c_0,c_1,\ldots,c_{n-1},c_n-2)$.
Now set $m_0=L$ and
\begin{eqnarray}
\lambda^{(\mu)}_i&=&{1\over 2}(m_{t_{\mu}+i}-m_{t_{\mu}+i+1}),\qquad
(1\le i<e_\mu,\quad 0\le\mu\le n);\label{LambdaDef1}\\
\lambda^{(\mu)}_{e_\mu}&=&{1\over 2}m_{t_{\mu+1}+1},\qquad
(0\le\mu<n);\label{LambdaDef2}\\
\lambda^{(n)}_{e_n}&=&{1\over 2}m_{t_{n+1}-1}.
\label{LambdaDef3}
\end{eqnarray}
In addition, set:
\begin{eqnarray}
w_0&=&L/2;\\
w_\mu&=&\lambda^{(\mu-1)}_{e_{\mu-1}},\qquad(1\le\mu\le n);\\
\lambda^{(\mu)}_0&=&w_{\mu-1}-wt(\lambda^{(\mu-1)}),\qquad(1\le\mu\le n)%
\label{LambdaDef6}
\end{eqnarray}
(It may be checked that (\ref{LambdaDef3}) and (\ref{LambdaDef6})
agree when $e_n=0$.)
Note that we don't require a $\lambda^{(0)}_0$ to be specified.
{}From these definitions, we obtain:
\begin{eqnarray}
m_{t_{\mu}+i}&=&2\left(w_\mu-\sum_{j=1}^{i-1}\lambda^{(\mu)}_j\right),
\qquad (1\le i\le e_\mu,\quad 0\le\mu\le n);\label{MintoLEq}\\
n_{t_{\mu}+i}&=&\lambda^{(\mu)}_{i-1}-\lambda^{(\mu)}_{i},\qquad
(1+\delta_{\mu,0}\le i\le e_\mu,\quad 0\le\mu\le n),\label{NintoLEq}
\end{eqnarray}
where the second expression arises on using Lemma \ref{CarLem}.
The restrictions $n_j\ge0$ are thus equivalent to
$\lambda^{(0)}_1\ge\lambda^{(0)}_2\ge\cdots\ge\lambda^{(0)}_{e_0}$
and
$\lambda^{(\mu)}_0\ge\lambda^{(\mu)}_1\ge\cdots\ge\lambda^{(\mu)}_{e_\mu}$
for $1\le\mu\le n$.
That $\lambda^{(\mu)}_{e_\mu}\ge0$ arises from (\ref{LambdaDef2})
and (\ref{LambdaDef3}).
Thus $\lambda^{(0)},\lambda^{(1)},\ldots,\lambda^{(n)}$
are indeed partitions with $\lambda^{(\mu)}_0\ge\lambda^{(\mu)}_1$
for $1\le\mu\le n$ ($\lambda^{(0)}_1$ is unbounded).
That $\lambda^{(\mu)}_0\ge0$ for $1\le\mu\le n$ implies, via
(\ref{LambdaDef6}), that $wt(\lambda^{(\mu)})\le w_\mu$
for $0\le\mu<n$.
Additionally, we obtain
$\sum_{i=1}^{e_n} \lambda^{(n)}_i={1\over 2}m_{t_n+1}=
\lambda^{(n-1)}_{e_{n-1}}=w_n$, so that $wt(\lambda^{(n)})=w_n$.

Thereupon, Theorem \ref{LFinThrm1} may be rewritten as follows
(cf.\ Theorem 3.8 of \cite{flw}):

\begin{theorem}\label{LFinThrm2}
Let $p$ and $p'$ be positive co-prime integers and $r_0$ and $s_0$
the smallest non-negative integers such that $\vert ps_0-p'r_0\vert=1$.
If $p'/p$ has continued fraction $(e_0,e_1,\ldots,e_{n-1},e_n+2)$,
then for even $L\ge0$:
\begin{eqnarray*}
\hbox to 0mm{$\displaystyle
\chi^{p,p'}_{s_0,s_0,s_0+1}(L)
$\hss}\\[0.5mm]
&&=\sum
q^{\sum_{\mu=0}^n\sum_{i=1}^{e_\mu} \lambda_{\mu,i}^2}
\prod_{\mu=0}^{n}\prod_{i=1+\delta_{\mu,0}}^{e_\mu}
\left[{2\left(w_\mu-\sum_{j=1}^i \lambda^{(\mu)}_{j}\right)
			  +\lambda^{(\mu)}_{i-1}+\lambda^{(\mu)}_{i}
\atop \lambda^{(\mu)}_{i-1}-\lambda^{(\mu)}_{i}}\right]_q
\!,
\end{eqnarray*}
where the sum is over all sequences
$\lambda^{(0)}$, $\lambda^{(1)},\ldots,\lambda^{(n)}$ of partitions for
which, for $0\le\mu\le n$,
the partition $\lambda^{(\mu)}=(\lambda^{(\mu)}_1,\lambda^{(\mu)}_2,
\ldots,\lambda^{(\mu)}_{e_\mu})$ satisfies
$\lambda^{(\mu)}_1\le\lambda^{(\mu)}_0$ ($0<\mu$)
and $wt(\lambda^{(\mu)})\le w_\mu$,
where we define $w_0=L/2$, $w_\mu=\lambda^{(\mu-1)}_{e_{\mu-1}}$
for $1\le\mu\le n$,
and $\lambda^{(\mu)}_0=w_{\mu-1}-wt(\lambda^{(\mu-1)})$
for $1\le\mu\le n$;
and additionally also satisfy $wt(\lambda^{(n)})=w_n$.
(In the above formula, $\lambda_{\mu,i}=\lambda^{(\mu)}_i$.)
\end{theorem}

\Proof On using (\ref{MintoLEq}), we readily find that:
$$
{1\over 4}\left( \boldm^T\boldC\boldm - L^2 \right)
= \sum_{\mu=0}^n\sum_{i=1}^{e_\mu} \lambda_{\mu,i}^2.
$$
The result then immediately follows from Theorem \ref{LFinThrm1},
using (\ref{MintoLEq}) and (\ref{NintoLEq}). 
\cqfd
\medskip

We are now able to take the $L\to\infinity$ limit.
The following theorem deals with the $p>2$ cases, for which
$n>0$ and $e_1>0$.

\begin{theorem}\label{LInfThrm}
Let $p$ and $p'$ be positive co-prime integers with $p'>p>2$,
and $r_0$ and $s_0$
the smallest non-negative integers such that $\vert ps_0-p'r_0\vert=1$.
If $p'/p$ has continued fraction $(e_0,e_1,\ldots,e_{n-1},e_n+2)$,
then:
\begin{eqnarray*}
\chi^{p,p'}_{r_0,s_0}
&&\!\!\!\!\!\!\!\!\!\!=\sum\left(
q^{\sum_{\mu=0}^n\sum_{i=1}^{e_\mu} \lambda_{\mu,i}^2}
{1\over(q)_{2\lambda^{(0)}_{e_0}}}
\prod_{i=2}^{e_0} {1\over(q)_{\lambda^{(0)}_{i-1}-\lambda^{(0)}_{i}}}\right.\\
&&\qquad\qquad\left.
\prod_{\mu=1}^{n}\prod_{i=1+\delta_{\mu,1}}^{e_\mu}
\left[{2\left(w_\mu-\sum_{j=1}^i \lambda^{(\mu)}_{j}\right)
			  +\lambda^{(\mu)}_{i-1}+\lambda^{(\mu)}_{i}
\atop \lambda^{(\mu)}_{i-1}-\lambda^{(\mu)}_{i}}\right]_q\; \right),
\end{eqnarray*}
where the sum is over all sequences
$\lambda^{(0)}$, $\lambda^{(1)},\ldots,\lambda^{(n)}$ of partitions for
which, for $0\le\mu\le n$,
the partition $\lambda^{(\mu)}=(\lambda^{(\mu)}_1,\lambda^{(\mu)}_2,
\ldots,\lambda^{(\mu)}_{e_\mu})$ satisfies
$\lambda^{(\mu)}_1\le\lambda^{(\mu)}_0$ ($1<\mu$)
and $wt(\lambda^{(\mu)})\le w_\mu$ ($0<\mu$),
where we define, $w_\mu=\lambda^{(\mu-1)}_{e_{\mu-1}}$
for $1\le\mu\le n$,
and $\lambda^{(\mu)}_0=w_{\mu-1}-wt(\lambda^{(\mu-1)})$
for $2\le\mu\le n$;
and additionally also satisfy $wt(\lambda^{(n)})=w_n$.
(In the above formula, $\lambda_{\mu,i}=\lambda^{(\mu)}_i$.)
\end{theorem}

\Proof By (\ref{ChiLimEq}),
$$
\chi^{p,p'}_{r_0,s_0} = \lim_{L\to\infinity}
\chi^{p,p'}_{s_0,s_0,s_0+1}(L),
$$
whereupon the result follows from Theorem \ref{LFinThrm2}
after noting that:
$$
\lim_{m\to\infinity}\left[{m+n\atop n}\right]_q
=\lim_{m\to\infinity}\left[{m+n\atop m}\right]_q
={1\over(q)_n}.
$$
\cqfd
\medskip

The case $p=2$ has $n=1$ and $e_1=0$. In this case, taking the
$L\to\infinity$ limit of the expression given by
Theorem \ref{LFinThrm2}, reproduces the summation expression
for $\chi^{2,2e_0+1}_{1,e_0}$ first given by Gordon \cite{gordon}
(see also \cite{andrews1,andrews-red-book,burge1,ab,bressoud}):

\begin{theorem}
$$
\chi^{2,2e_0+1}_{1,e_0}
=\sum_{\lambda}
{q^{\lambda_1^2+\lambda_2^2+\cdots+\lambda_{e_0-1}^2}\over
(q)_{\lambda_1-\lambda_2}(q)_{\lambda_2-\lambda_3}\cdots
(q)_{\lambda_{e_0-2}-\lambda_{e_0-1}}
(q)_{\lambda_{e_0-1}} },
$$
where the sum is over all
partitions $\lambda=(\lambda_1,\lambda_2,\ldots,\lambda_{e_0-1})$.
\end{theorem}

\Proof With $p=2$ and $p'=2e_0+1$, the continued fraction of
$p'/p$ is $(e_0,2)$, so that $e_1=0$.
Then Theorem \ref{LFinThrm2} requires a sum over all pairs
of partitions $\lambda^{(0)}$ and $\lambda^{(1)}$ which have
$e_0$ and $e_1$ parts respectively. Since it is required that
$wt(\lambda^{(1)})=w_1$ and $w_1=\lambda^{(0)}_{e_0}$,
it is necessary that $\lambda^{(0)}_{e_0}=0$.
The only other constraint is $wt(\lambda^{(0)})\le L/2$.
The result now follows exactly as in the proof of Theorem
\ref{LInfThrm}.
\cqfd
\medskip

In the $p=1$ case, the sum of Theorem \ref{LFinThrm2} is over
all partitions $\lambda^{(0)}$ such that $wt(\lambda^{(0)})=L/2$.
Thus, as $L\to\infinity$, the exponent of $q$ is unbounded and
so $\lim_{L\to\infinity}\chi^{1,p'}_{1,1,2}(L)$ does not exist.
This is to be expected since as $L\to\infinity$,
every path in ${\P}^{1,p'}_{1,1,2}(L)$ has an unbounded
number of scoring vertices.

\begin{appendix}

\section{Bijection between paths and partitions}

In this Appendix, we describe a natural weight preserving bijection 
between the Forrester-Baxter paths $\P^{p, p'}_{a, b, c}(L)$ \cite{fb}, 
and partitions with prescribed hook differences \cite{abbbfv}. That 
such a bijection exists was anticipated in \cite{abbbfv}, where 
generating functions for the two are shown to be equal up to 
a normalisation.

\subsection{Path recurrence relations}

For $0<h<p'$, define the functions $r(h)=\lfloor ph/p'\rfloor$ and 
$\hat r(h)=\lfloor (p'-p) h / p'\rfloor$ and note that 
$r(h)+\hat r(h)=h-1$. From the definitions (\ref{BwtDef}) and 
(\ref{BgenDef}), it is easy to see that the path generating functions 
$\phi^{p, p'}_{a,b,c}(L)$ satisfy the recurrences (\cite{fb}):
\begin{eqnarray}
\phi^{p, p'}_{a,b,b+1}(L)
&=&q^{L\hat r(b+1)} \phi^{p, p'}_{a,b+1,b}(L-1)
+q^{L/2} \phi^{p, p'}_{a,b-1,b}(L-1),
\label{Brec1}\\[0.5mm]
\phi^{p, p'}_{a,b,b-1}(L)
&=&q^{L/2} \phi^{p, p'}_{a,b+1,b}(L-1)
+q^{-L\hat r(b-1)} \phi^{p, p'}_{a,b-1,b}(L-1);
\label{Brec2}
\end{eqnarray}
the boundary conditions:
\begin{eqnarray}
\phi^{p, p'}_{a,b-1,b}(L)=0&&\hbox{if $b=1$,}\\[0.5mm]
\phi^{p, p'}_{a,b+1,b}(L)=0&&\hbox{if $b=p'-1$;}
\end{eqnarray}
and initial conditions:
\begin{equation}
\phi^{p, p'}_{a,b,b+1}(0)=\phi^{p, p'}_{a,b,b-1}(0)=\delta_{a,b}.
\end{equation}
These five properties uniquely determine $\phi^{p, p'}_{a,b,c}(L)$
in all cases.  Now define
\begin{equation}\label{ChiLDef}
\chi^{p,p'}_{a,b,b\pm1}(L)
=q^{-\hat r(c)(a-b\pm L)/2-(a-b)(a-c)/4} \phi^{p,p'}_{a,b,b\pm1}(L),
\end{equation}
where $c=b\pm1$.
In translating the above recurrences to this new function,
it is appropriate to treat the two cases of
$\hat r(b)=\hat r(b\pm1)$ and $\hat r(b)\ne \hat r(b\pm1)$ separately.
In the $\hat r(b)=\hat r(b\pm1)$ case, we find that (\ref{Brec1}) and
(\ref{Brec2}) become
\begin{equation}\label{Rec1}
\chi^{p,p'}_{a,b,b+1}(L)
=\chi^{p,p'}_{a,b+1,b}(L-1)
+q^{(L+a-b)/2}\chi^{p,p'}_{a,b-1,b}(L-1),
\end{equation}
and
\begin{equation}\label{Rec2}
\chi^{p,p'}_{a,b,b-1}(L)
=q^{(L-a+b)/2}\chi^{p,p'}_{a,b+1,b}(L-1)
+\chi^{p,p'}_{a,b-1,b}(L-1),
\end{equation}
respectively.
In the other case, where necessarily $\hat r(b\pm1)=\hat r(b)\pm1$, we have
\begin{equation}\label{Rec3}
\chi^{p,p'}_{a,b,b+1}(L)
=q^{(L-a+b)/2}\chi^{p,p'}_{a,b+1,b}(L-1)
+\chi^{p,p'}_{a,b-1,b}(L-1),
\end{equation}
and
\begin{equation}\label{Rec4}
\chi^{p,p'}_{a,b,b-1}(L)
=\chi^{p,p'}_{a,b+1,b}(L-1)
+q^{(L+a-b)/2}\chi^{p,p'}_{a,b-1,b}(L-1).
\end{equation}
The boundary and initial conditions are similar to those above.
Namely:
\begin{eqnarray}
\chi^{p,p'}_{a,b-1,b}(L)=0&&\hbox{if $b=1$,}\\[0.5mm]
\chi^{p,p'}_{a,b+1,b}(L)=0&&\hbox{if $b=p'-1$,}
\end{eqnarray}
and
\begin{equation}
\chi^{p,p'}_{a,b,b+1}(0) = \chi^{p,p'}_{a,b,b-1}(0) = \delta_{a,b}.
\end{equation}

\subsection{Partitions with prescribed hook differences}\label{PartSec}

A partition $\mu=(\mu_1,\linebreak[3]\mu_2,\ldots,\mu_M)$ is a sequence of $M$
integer parts $\mu_1,\mu_2,\ldots,\mu_M,$ satisfying
$\mu_1\ge\mu_2\ge\cdots\ge\mu_M>0$.
It is to be understood that $\mu_i=0$ for $i>M$.
The weight $\wt(\mu)$ of $\mu$ is given by $\wt(\mu)=\sum_{i=1}^M\mu_i$.
The partition $\mu$ is often depicted by its {\it Young diagram}
(also called Ferrars graph),
$F^\mu$ which comprises $M$ left-adjusted
rows, the $i$th row of which (reading down) consists of $\mu_i$
cells \cite{andrews-red-book}.
The coordinate $(i,j)$ of a cell is obtained by setting $i$
and $j$ to be respectively, the row and column (reading from the left)
in which the cell resides.
The $k$th diagonal of $F^\mu$ comprises all those cells of $F^\mu$
with coordinates $(i,j)$ which satisfy $i-j=k$.

The partition $\mu^\prime$, conjugate to $\mu$, is obtained
by setting $\mu^\prime_j$ to be the number of cells in the
$j$th column of $F^\mu$.
The hook difference at the cell with coordinate $(i,j)$
is then defined to be $\mu_i-\mu^\prime_j$.
As an example, filling each cell of $F^{(5,4,3,1)}$ with its
hook difference, yields:
$$
\youngd{
\multispan{11}\hrulefill\cr
&1&&\mathbf2&&2&&3&&4&\cr
\multispan{11}\hrulefill\cr
&0&&1&&\mathbf1&&2&\cr
\multispan{9}\hrulefill\cr
&-1&&0&&0&\cr
\multispan{7}\hrulefill\cr
&-3&\cr
\multispan{3}\hrulefill\cr}\;.
$$
The bold entries are those on diagonal $-1$.
In what follows, we will be especially interested in the
hook differences on certain diagonals.

Let $K,i,N,M,\alpha,\beta$ be non-negative integers for
which $1\le i\le K/2$, $\alpha+\beta<K$ and
$\beta-i\le N-M\le K-\alpha-i$.
In \cite{abbbfv}, $D_{K,i}(N,M;\alpha,\beta)$ is defined to be
the generating function for partitions $\mu$
into at most $M$ parts, each not exceeding $N$ such that
the hook differences on diagonal $1-\beta$ are at least
$\beta-i+1$, and on diagonal $\alpha-1$ are at most $K-i-\alpha-1$.
In addition, if $\alpha=0$, the restriction that $\mu_{N-L+i+1}>0$
is also imposed; and if $\beta=0$, the restriction that
$\mu_1>M-i$ is also imposed.
It may be shown (see \cite[page 346]{abbbfv}) that:
\begin{eqnarray}
&&
D_{K,i}(N,M;\alpha,\beta)
= 
D_{K,i}(N,M\!-\!1;\alpha,\beta)
+q^M D_{K,i}(N\!-\!1,M;\alpha\!+\!1,\beta\!-\!1),
\label{ParRec1}\\[0.5mm]
&&
D_{K,i}(N,M;\alpha,\beta)
= 
D_{K,i}(N\!-\!1,M;\alpha,\beta)
+q^N D_{K,i}(N,M\!-\!1;\alpha\!-\!1,\beta\!+\!1),
\label{ParRec2}
\end{eqnarray}
\begin{eqnarray*}
D_{K,i}(M+K-i,M;0,\beta)&=&0,\\[0.5mm]
D_{K,i}(M-i,M;\alpha,0)&=&0,\\[0.5mm]
D_{K,i}(0,0;\alpha,\beta)&=&1.
\end{eqnarray*}
Now if we define
\begin{displaymath}
\hat\chi^{p,p'}_{a,b,c}(L)=
D_{p',a}\left({L-a+b\over 2},{L+a-b\over 2};p-r,r\right),
\end{displaymath}
with $r=\lfloor pc/p'\rfloor+(b-c+1)/2$.
we find that $\hat\chi^{p,p'}_{a,b,c}(L)$ satisfies precisely
the same recursion, boundary and initial conditions as
$\chi^{p,p'}_{a,b,c}(L)$. Since these conditions determine
a function uniquely, we conclude that
$\hat\chi^{p,p'}_{a,b,c}(L)=\chi^{p,p'}_{a,b,c}(L)$.
That is, we obtain:
\begin{equation}\label{Compare}
\chi^{p,p'}_{a,b,c}(L)=
D_{p',a}\left({L-a+b\over 2},{L+a-b\over 2};p-r,r\right).
\end{equation}

In \cite{abbbfv}, the above recurrences for
$D_{K,i}(N,M;\alpha,\beta)$ are solved to yield
\begin{eqnarray*}
D_{K,i}(N,M;\alpha,\beta)&=&
\sum_{\lambda=-\infinity}^\infinity
q^{\lambda(K\lambda-i)(\alpha+\beta)+K\beta\lambda}
\left[ {N+M\atop M-K\lambda} \right]_q\\[0.5mm]
&&\qquad\quad
-\sum_{\lambda=-\infinity}^\infinity
q^{\lambda(K\lambda+i)(\alpha+\beta)+K\beta\lambda+\beta i}
\left[ {N+M\atop M-K\lambda-i} \right]_q.
\end{eqnarray*}
Using this and (\ref{Compare}), we are led to the expression
for $\chi^{p,p'}_{a,b,c}(L)$ given in (\ref{FinRochaEq}).

\subsection{The bijection}

In order to describe the bijection between the paths and the
partitions described above, we recapitulate the derivation
of (\ref{ParRec1}) and (\ref{ParRec2}) given in \cite{abbbfv}.
Consider one of the partitions enumerated by $D_{K,i}(N,M;\alpha,\beta)$.
If the $M$th part is zero then the partition
also appears in $D_{K,i}(N,M-1;\alpha,\beta)$.
On the other hand, if the $M$th part is at least 1,
then the partition obtained by decreasing each part by 1
appears in $D_{K,i}(N-1,M;\alpha+1,\beta-1)$
(the removal of the first part changes the diagonals on which the
hook differences are considered --- the changes in $\alpha$ and
$\beta$ reflect this)
Thus, (\ref{ParRec1}) results.

Now, under the identification (\ref{Compare}), expression
(\ref{Rec1}) is equivalent to (\ref{ParRec1}).
The first term of (\ref{Rec1}) corresponds to a peak-down
vertex at the $L$th position and the second term to a straight-up
vertex at the $L$th position.
This indicates that if $\hat r(b)=\hat r(b+1)$
(so that $r(b)\ne r(b+1)$ and hence the second edge of the vertex is
in an odd band), then a straight-up vertex
corresponds to being able to remove the first column of
length $M=(L+a-b)/2$ from the partition. On the other hand,
if $\hat r(b)=\hat r(b+1)$, then a peak-down vertex corresponds to
leaving the partition unchanged, but reducing the constraining
parameter $M$ by 1.

When $\hat r(b)\ne \hat r(b-1)$, the analysis of (\ref{Rec4}) is similar.
In this case, a peak-up vertex
corresponds to being able to remove the first column of
length $M=(L+a-b)/2$ from the partition, and
a peak-down vertex corresponds to
leaving the partition unchanged, but again reducing the constraining
parameter $M$ by 1. We thus obtain the first and fourth entries
in the first column of the following table.

Again consider one of the partitions enumerated by
$D_{K,i}(N,M;\alpha,\beta)$.
If the first part is less than $N$ then the partition also appears
in $D_{K,i}(N-1,M;\alpha,\beta)$.
On the other hand, if the first part is exactly $N$, then
the partition obtained by removing the first part
appears in $D_{K,i}(N,M-1;\alpha-1,\beta+1)$.
Thus, (\ref{ParRec2}) results.

Under the identification (\ref{Compare}), expression
(\ref{Rec2}) is equivalent to (\ref{ParRec2}).
The first and second terms of (\ref{Rec2}) correspond respectively
to a straight-down vertex and a peak-up vertex at the $L$th position.
Thus, when $\hat r(b)=\hat r(b-1)$, a straight-down vertex corresponds
to being able to remove the first part of length $N=(L-a+b)/2$
from the partition, and a peak-up vertex corresponds to
leaving the partition unchanged, but reducing the constraining
parameter $N$ by 1.

For $\hat r(b)\ne \hat r(b+1)$, a similar analysis of (\ref{Rec3})
shows that a peak-down vertex corresponds to being able to
remove the first part of length $N=(L-a+b)/2$ from the
partition, and a straight-up vertex corresponds to
reducing the parameter $N$ by 1.

By recursively applying the above rules, through a process
of successive row and column removal, we eventually arrive
at the empty partition corresponding to a path of length 0.
By applying the procedure in reverse, traversing the path from
left to right and keeping track of
the values $N$ and $M$ as we proceed, we can build the
required partition up from the empty partition.
For this, the above description yields the moves given in Table 1.

\begin{table}[htb]
\begin{center}
\begin{tabular}{|c|c|c|c|c|c|}
\hline
Vertex&
      &
Move  &
Vertex&
      &
Move  \\
\hline
\hline\wombat
\epsfbox{v1.eps}&
\raisebox{5mm}{\parbox{25mm}{ Add $M$ column \newline \& increment $N$}}&
\raisebox{3.5mm}{\epsfbox{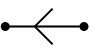}}&
\epsfbox{v5.eps}&
\raisebox{5mm}{\parbox{25mm}{ Increment $N$}}&
\raisebox{3.5mm}{\epsfbox{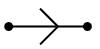}}\\
\hline\wombat
\epsfbox{v2.eps}&
\raisebox{5mm}{\parbox{25mm}{ Add $N$ row \newline \& increment $M$}}&
\epsfbox{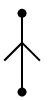}&
\epsfbox{v6.eps}&
\raisebox{5mm}{\parbox{25mm}{ Increment $M$}}&
\epsfbox{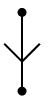}\\
\hline\wombat
\epsfbox{v3.eps}&
\raisebox{5mm}{\parbox{25mm}{ Increment $N$}}&
\raisebox{3.5mm}{\epsfbox{arrow_rt.eps}}&
\epsfbox{v7.eps}&
\raisebox{5mm}{\parbox{25mm}{ Add $M$ column \newline \& increment $N$}}&
\raisebox{3.5mm}{\epsfbox{arrow_lf.eps}}\\
\hline\wombat
\epsfbox{v4.eps}&
\raisebox{5mm}{\parbox{25mm}{ Increment $M$}}&
\epsfbox{arrow_dn.eps}&
\epsfbox{v8.eps}&
\raisebox{5mm}{\parbox{25mm}{ Add $N$ row \newline \& increment $M$}}&
\epsfbox{arrow_up.eps}\\
\hline
\end{tabular}
\end{center}
\vskip2mm
\centerline{Table 1.}
\end{table}

Alongside each description, we give a line segment and an arrow
which encapsulates the description. This provides a handier
means of building the partition. In fact, it describes a construction
of the partition's profile (outline) as follows.
Let the length zero path correspond to a single dot. Then,
scanning the path from left to right, for each vertex,
append the appropriate line segment to
the profile so that the arrow points away from the
initial dot. After the $L$th vertex has been considered,
this profile is abutted into the axes in the 4th quadrant to
produce the Young diagram of the required partition.

To illustrate this procedure, consider the following path in 
$\P^{3,8}_{4,3,2}(15)$:

\medskip
\begin{center}
\epsfig{file=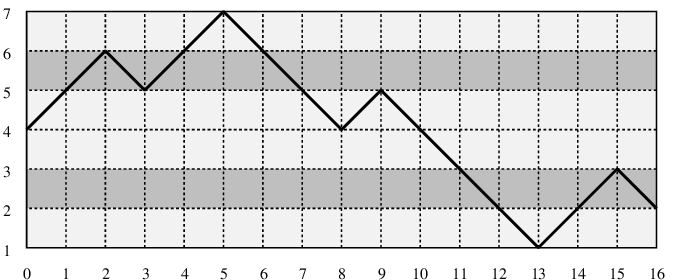}
\end{center}
\medskip
 
\noindent
The above procedure then produces the following Young diagram
(here the starting point is the unfilled circle).

\medskip
\begin{center}
\epsfig{file=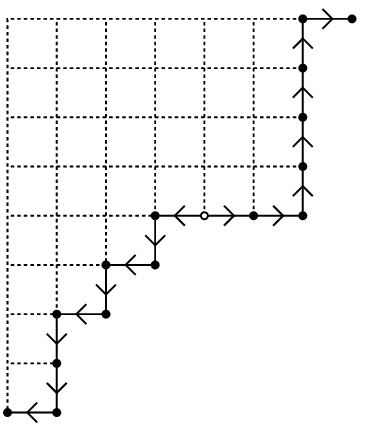}
\end{center}
\medskip

\noindent
Thus the required partition is $(6,6,6,6,3,2,1,1)$.

Using the above table, we see that the values of $M$ and $N$
at any vertex count the number of SE segments and NE segments
respectively prior to that vertex. Thus in terms of the
$xy$-coordinate system defined in Section \ref{PresSec}, we simply
have $x=M$ and $y=N$. The definition (\ref{WtDef}) then gives
the weight (sum of parts) of the corresponding partition.

That the map from path to partition is in fact a bijection
follows because given a partition and values of $a,b,c$ and $L$,
we obtain $N=(L-a+b)/2$ and $M=(L+a-b)/2$ at the $L$th vertex.
Then one of the two extremal line segments of the partition's
profile (extended to be of width $N$ and height $M$) corresponds
to the $L$th vertex. The right edge of the vertex is determined
by $b$ and $c$, whereupon Table~2 shows that only
one left edge can occur. After removing the corresponding line
segment from the partition's profile, this process is
recursively repeated to produce a unique complete path.

We note that in the case where $p'=p+1$, which implies that all
bands are odd, and where $a=b=1$, this bijection reduces to
that given in \cite{fw}.
In the case $p=2$, it reduces to that given in
\cite{burge2,ab,bressoud}.

\end{appendix}

\subsection*{Acknowledgements}

We wish to thank Iain Aitchison, Alex Feingold, Dominique Foata, 
Christian Krattenthaler, Tetsuji Miwa, Dennis Stanton and Ole 
Warnaar for discussions, and for their interest in this work. We 
also wish to thank the Australian Research Council for financial 
support.

\end{document}